\documentclass[12pt, oneside]{amsart}
\usepackage{geometry}                
\usepackage{graphicx}
\usepackage{amssymb}

\usepackage{amsthm}
\usepackage{amscd}
\usepackage[all]{xy}
\usepackage{pinlabel}
\usepackage[abs]{overpic}
\usepackage{comment}
\usepackage{subcaption}
\usepackage[dvipsnames]{xcolor} 

\newcommand{\e}{\varepsilon}
\newcommand{\omu}{\overline{\mu}}

\theoremstyle{theorem}
\newtheorem{thm}{Theorem}[section]
\newtheorem{prop}[thm]{Proposition}
\newtheorem{cor}[thm]{Corollary}
\newtheorem{lem}[thm]{Lemma}

\theoremstyle{definition}
\theoremstyle{remark}
\newtheorem{remark}[thm]{Remark}
\newtheorem*{acknowledgments}{Acknowledgments}

\begin{document}


\title[Milnor Invariants]{Milnor Invariants \\ 
\footnotesize{From classical links to surface-links, and beyond}}

\author{Akira Yasuhara}
\address{Faculty of Commerce, Waseda University, 1-6-1 Nishi-Waseda,
  Shinjuku-ku, Tokyo 169-8050, Japan}
	 \email{yasuhara@waseda.jp}

\maketitle
\section*{Introduction}\label{sec1}
This is an English translation of the expository article written by the author in Japanese 
for publication in {\em Sugaku}.

Milnor invariants for (classical) links were defined by J. Milnor \cite{M1},\cite{M2} in the 1950's.
The study of Milnor invariants has a long history, and counts numerous research results. 
It is impossible (at least for the author) to cover all of them, so in this article 
the author only try to explain Milnor invariants from the viewpoint of his research.
In \cite{JB-Dic}, J.B. Meilhan explained Milnor invariants 
from a different perspective, and gave a concise summary 
of the topics not covered here.

This is not a research paper, so we shall sometimes sacrifice precision and give rough 
explanations so that a broader audience could follow the outline of this article.
On the other hand, we also cover topics that are oriented towards readers who are familiar with knot theory, 
so we may use technical terms without explanation.
For explanations, we sometimes use different terminologies and notations from those in the references.
In addition, we occasionally  simplify statements of results to fit this article.

The content of this article is as follows:
In chapter~1, we treat Milnor invariants for {\em classical} links.
A classical link is a mathematical model of a \lq closed strings in space\rq. 
Hence it is closely related to real-world object.
On the other hand, by ignoring the real world and focusing only on the structure of the link, 
we obtain something called a {\em welded link}. 
In chapter~2, we will discuss Milnor invariants for welded links.
We also present an algorithm for computing Milnor invariants.
The relationship between \lq classical links\rq\ and  \lq welded links\rq\  is 
in some sense similar to the relationship between \lq real numbers' and \lq complex numbers\rq\, 
and in fact, classical links are  embedded in welded links.
Extending the study of classical links to welded links may lead to new discoveries.
In chapter~3 we present results giving geometric characterizations of 
Milnor invariants of welded links; 
these results also hold for classical links, but the results cannot be derived 
from observations of classical links.
In chapter ~4, for surfaces in 4-space, i.e., for {\em surface-links}, 
we introduce  {\em $2$-dimensional cut-diagrams}, which can be regarded as welded surface-links.
We define Milnor invariants of $2$-dimensional cut-diagrams, and via this definition, 
we also define Milnor invariants of surface-links.   
In addition, we present an algorithm for computing these Milnor invariants.
The definition of a 2-dimensional cut-diagram can be extended to that of 
an $m$-dimensional cut-diagram $(m\geq 2)$.
In chapter~5, we define Milnor invariants for $m$-dimensional cut-diagrams, and 
likewise for $m$-dimensional links. 

\begin{acknowledgments}
The author thanks Jean-Baptiste Meilhan for useful comments on a draft version of this article.  
\end{acknowledgments}

\section{Milnor invariants for classical links}

In this chapter, we give a quick overview of Milnor invariants for classical links.

\medskip
\subsection{Links and string links}

Let $n$ be a positive integer. An {\em $n$-component link} is a union of  
$n$ simple closed curves in 3-space.  In particular, a 1-component link is called a  {\em knot}. 
For example, in Figure~\ref{link}, $K_1\cup K_2$ is a $2$-component link. 
A link is {\em trivial} if it bounds a disjoint union of disks. 
Two links are {\em equivalent} if there is a \lq continuous deformation\rq\  
(more precisely, {\em an ambient isotopy}) between them.


 Let $D^2$ be the unit disk and $D^1$ the diameter on the $x$-axis with same orientation 
 as the $x$-axis. 
 An {\em $n$-string link} is a union of  
 $n$ simple curves in the cylinder $D^2\times [0,1]$ 
 such that the $i$th component runs from $p_i\times\{0\}$ to $p_i\times\{1\}$ for each $i=1,...,n$,
 where $p_1,...,p_n$ are points on $D^1$ that are arranged in order along the orientation of $D^1$.
For example, in Figure~\ref{string}, $K_1\cup K_2\cup K_3$ is a $3$-string link. 

Two strings links are {\em equivalent} if there is a continuous deformation  between them 
fixing the boundary of  $D^2\times[0,1]$. 
An $n$-string link is {\em trivial} if it is equivalent to the $n$-string link 
$\bigcup_{i=1}^n(p_i\times[0,1])$.

\medskip
\subsection{Peripheral system and Milnor invariants for links}
\begin{figure}[!h]
\begin{center}
\includegraphics[width=.3\linewidth]{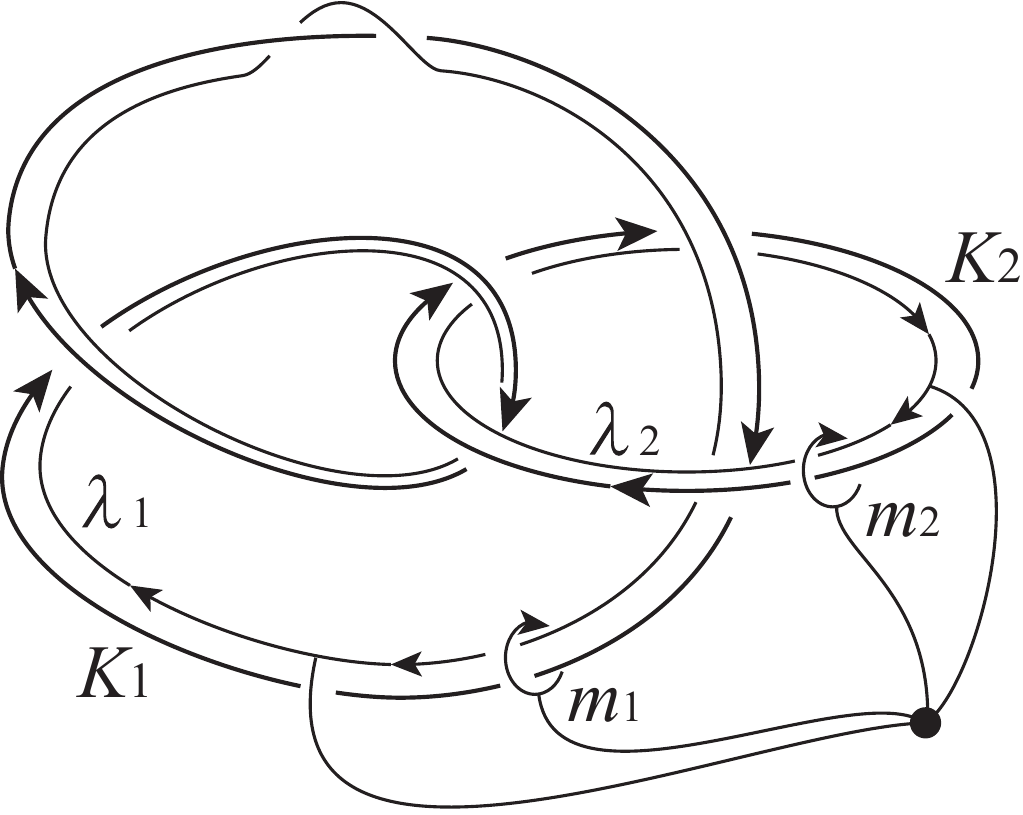}
   \vspace{.5em}
\caption{Meridian and longitude of a link}\label{link}
\end{center}
\end{figure}

For an $n$-component link $L=K_1\cup\cdots\cup K_n$, let $G(L)$ be the fundamental group 
of the complement of $L$. 
As illustrated in Figure~\ref{link}, for each component $K_i$, we choose a pair of 
elements $m_i$ and $\lambda_i$ in $G(L)$, that are called an $i$th {\em meridian} and an $i$th 
{\em longitude} respectively,\footnote{While we skip the detailed definition, in this article it is enough to know  that meridians and longitudes are special elements in $G(L)$ that depend on the link. } and 
we call the pair $(G(L),\{m_i,\lambda_i\}_i)$ a   {\em peripheral system}  
 of $L$.  
It is known that the peripheral systems up to certain equivalence give the classification of links \cite[Theorem 6.1.7]{Kaw}. 
Hence investigating the equivalence classes of peripheral systems is nothing more than 
investigating  the equivalence classes of links. 
As we will see below, Milnor invariants are invariants derived from peripheral systems. 

Since $G(L)$ is non-commutative, classifying peripheral systems is as difficult as classifying links.
Therefore we consider the natural projection 
\[\rho_2:G(L)\longrightarrow N_2(L)=G(L)/\Gamma_2G(L),\]
where $\Gamma_2G(L)=[G(L),G(L)]$ is the commutator subgroup\footnote{$[X,Y]$ is the subgroup generated 
by the set $\{xy^{-1}x^{-1}y~|~x\in X,y\in Y\}$}
of $G(L)$. 
We remark that $N_2(L)$ is a free abelian group generated by $\rho_2(m_1),...,\rho_2(m_n)$. 
Hence each $\rho_2(\lambda_k)~(k=1,...,n)$ can be written as 
\[\rho_2(\lambda_k)=\sum_{i\in\{1,...,n\}\setminus\{ k\}} \mu_L(ik)\rho_2(m_i).\]
The integer coefficient $\mu_L(ik)$ of $\rho_2(m_i)$ is called a {\em (length-2) 
Milnor invariant} of $L$. 

More generally, for an integer $q\geq 2$ we consider the natural projection
\[\rho_q:G(L)\longrightarrow N_q(L)=G(L)/\Gamma_qG(L),\]
where for a group $G$, we define $\Gamma_1G=G$ and $\Gamma_qG=[\Gamma_{q-1}G,G]$. 
It is known that $N_q(L)$ is a nilpotent group generated by $\rho_q(m_1),...,\rho_q(m_n)$, 
and it is called the {\em $q$th nilpotent group} of $L$ (or the {\em $q$th nilpotent quotient} of $G(L)$). 
Taking the {\em Magnus expansion} $E(\rho_q(\lambda_k))$, we have 
an integer $\mu_L(I)$ with respect to each sequence $I=i_1,...,i_sk~(s<q)$ with length $s+1$. 
Here the Magnus expansion $E(\rho_q(\lambda_k))$ is a formal power series in non-commutative 
variables $X_1,...,X_n$ given by 
\[E(\rho_q(m_i))=1+X_i,~E(\rho_q(m_i^{-1}))=1-X_i+X_i^2-X_i^3+\cdots \hspace{1em}(i=1,2,\ldots,n).\]
We denote the integer coefficient of $X_{i_1}\cdots X_{i_s}$ in $E(\rho_q(\lambda_k))$ 
by $\mu_L(i_1...i_sk)$.  And we define $\mu_L(k)=0$.

The Magnus expansion has the following property.

\medskip
\begin{prop}\label{Magnus}(\cite{MKS}, \cite{Fenn})
Let $F$ be the free group generated by $\alpha_1,...,\alpha_n$, and let $E$ be the Magnus expansion defined by 
 $E(\alpha_i)=1+X_i,~E(\alpha_i^{-1})=1-X_i+X_i^2-X_i^3+\cdots$ for each $i$. Then the followings hold.
\begin{enumerate}
\item $E:F\longrightarrow E(F)$ is a bijection, and for any elements $g,h\in F$, $E(gh)=E(g)E(h)$.
\item For an element $g\in F$, $g\in\Gamma_qF$ if and only if the minimal degree of $E(g)-1$ is at least $q$. 
\end{enumerate}
\end{prop}

\medskip
For a link $L$, the numbers $\mu_L(I)$ are determined by a peripheral system  $(G(L),\{m_i,\lambda_i\}_i)$. 
But $(G(L),\{m_i,\lambda_i\}_i)$ is not uniquely determined by $L$, and $\mu_L(I)$ 
(of length at least 3) are not invariants for $L$. So, for a sequence $I$ with length at most $q$, we define 
\[\Delta_L(I):=\gcd
\left\{\mu_L(J)~\vline~
\begin{array}{l}
J\textrm{ is a sequence obtained from $I$ by deleting at least one}\\ 
\textrm{index and permuting the resulting sequence cyclicly}
\end{array}
\right\},\]
and take the residue class $\overline{\mu}_L(I)$ of $\mu_L(I)$ modulo $\Delta_L(I)$, which 
is an invariant for $L$ and is called a \emph{Milnor $\overline{\mu}$-invariant}.\footnote{A length-$2$ 
Milnor invariant $\mu_L(ik)$ is equal to the {\em linking number} $\mathrm{lk}(K_i,K_k)$
 between $K_i$ and $K_k$. Hence Milnor invariants are regarded as a generalization of linking number. 
V.G. Turaev \cite{Tur} and R. Porter \cite{Por} characterized 
Milnor invariants by means of {\em Massey product}. 
Moreover, T. Cochran \cite{Co} gave a geometric characterization 
for Milnor invariants as linking numbers of intersections 
between {\em Seifert surfaces} for links, which are closed related to Massey product. 
K. Murasugi \cite{Mur} showed that Milnor invariants are given by linking numbers of 
{\em branched covering space}.}
We note that if the length of $I$ is $2$, then $\Delta_L(I)=0$, and hence $\overline{\mu}_L(I)=\mu_L(I)$.
The length of the sequence $I$ is called \emph{the length of Milnor invariant} $\overline{\mu}_L(I)$.
It seems that Milnor invariants depend on $q$, but 
it is known that Milnor invariants of length at most $q$ are equal to those 
obtained from $\rho_{q+1}$. 
Therefore for any sequence  $I$ with length at most $q$, 
Milnor invariant $\overline{\mu}_L(I)$ is independent of $q$. 

\begin{figure}[!h]
\begin{center}
\includegraphics[width=.25\linewidth]{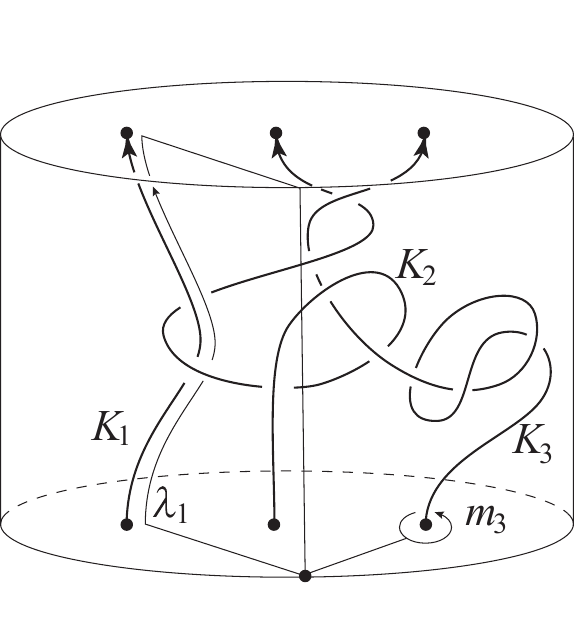}
   \vspace{.5em}
\caption{Meridian and longitude for a string link}\label{string}
\end{center}
\end{figure}

As in the case of links, we obtain the integer $\mu_L(I)$ 
for a string link $L$. 
In this case, the meridians and the longitudes are uniquely chosen as illustrated in Figure~\ref{string}. 
Hence for all $I$, $\mu_L(I)$ are invariants for $L$, and they are  called
{\em Milnor $\mu$-invariants}. 

\medskip
\subsection{Automorphisms of nilpotent groups and Milnor invariants}

Milnor invariants for string links are defined by N. Habegger and X.S. Lin~\cite{HL}. 
Their definition is different from that in the previous section. 
In this section, we present the idea of the definition by Habegger and Lin. 

Set $D_{\varepsilon}=D^2\times\{\varepsilon\}~(\varepsilon\in\{0,1\})$. 
For an $n$-string link $L\subset D^2\times[0,1]$, by Stallings Theorem~\cite[Theorem 5.1]{St}, 
the inclusion map 
\[D_{\varepsilon}\setminus L\hookrightarrow 
D^2\times[0,1]\setminus L\]
induces the isomorphism
\[\varphi_q^{\varepsilon}:\pi_1(D_{\varepsilon }\setminus L)/
\Gamma_q\pi_1(D_{\varepsilon}\setminus L)
\longrightarrow N_q(L)\]
for each $q$.
Hence we have an isomorphism  
\[\varphi_q:=(\varphi_q^1)^{-1}\circ \varphi_q^0:\pi_1(D_{0}\setminus L)/
\Gamma_q\pi_1(D_{0}\setminus L)
\longrightarrow \pi_1(D_{1}\setminus L)/
\Gamma_q\pi_1(D_{1}\setminus L).
\]
Since  
$\pi_1(D_{\varepsilon}\setminus L)(=\pi_1(D^2\setminus\{p_1,...,p_n\}))$ is the free group $F$ generated by 
the meridians $m_1,...,m_n$,  
we have the automorphism 
\[\varphi_q:F/
\Gamma_q F \longrightarrow F/
\Gamma_q F.\]
It follows from the definition of $\varphi_q$ that $\varphi_q(\rho_q(m_i))=\rho_q(\lambda_i^{-1}m_i\lambda_i)$. 
This means that $\varphi_q$ \lq contains\rq\  the information of Milnor $\mu$-invariants for $L$.  
In fact, $\varphi_q$ can be regarded as Milnor invariants of length at most $q-1$, see Proposition~\ref{meilhan-yasu-lemma}.

\medskip
\subsection{Properties of Milnor invariants} 

We summarize the properties of Milnor invariants 
that are deeply related to the topics in this article.

\begin{enumerate}
\item (J. Milnor \cite{M2}, J. Stallings \cite{St}, A. J. Casson \cite{Casson})~
Milnor invariants are \emph{isotopy} invariants for links \cite{M2}, and moreover   
\emph{link concordance} invariants \cite{St}, \cite{Casson}. 

\item (K. Habiro \cite{H})~Milnor invariants with length at most $k$ are 
\emph{$C_k$-equivalence} invariants, where 
the $C_k$-equivalence is an equivalence relation generated by local moves, 
\emph{$C_k$-moves}, defined by Habiro \cite{H}. 

\item (J. Milnor \cite{M1}, N. Habegger and X.S. Lin \cite{HL})~
For any sequence $I$ consisting of non repeated indices,  
$\overline{\mu}(I)$ is \emph{link-homotopy} invariant \cite{M1}, where 
link-homotopy is an equivalence relation generated by changing
\emph{crossings} between strands of the same component. 
It is known that 
$L$ is link-homotopic to a trivial link if and only if $\overline{\mu}_L(I)=0$ 
for any non-repeated sequence $I$ \cite{M1}, and that the $\mu(I)$'s classify string links up to 
link-homotopy \cite{HL}.\footnote{This follows from \cite[Theorem 1.7]{HL}, 
which is not stated in terms of Milnor invariants.}

\end{enumerate}

\section{Diagrams and Milnor invariants}

In this chapter, we define Milnor invariants for welded links, which are
defined as a generalization of (diagrams of) classical links,  and we also 
give an algorithm for computing these Milnor invariants. 

\medskip
\subsection{Link diagrams}
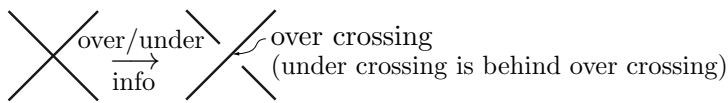
\begin{figure}[!h]
  \begin{center}
    \begin{overpic}[width=3.6cm]{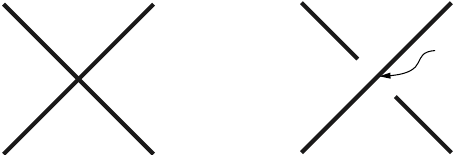}
      \put(27,20){\footnotesize over/under}
      \put(39,12){$\longrightarrow$}
      \put(40,4){\footnotesize info}
      \put(100,21){\small over crossing}
      \put(100,10){\footnotesize (under crossing is behind over crossing)}
    \end{overpic}
  \end{center}
  \caption{Crossing}
  \label{crossing}
\end{figure}

A link $L$ is an object in 3-space, but it can actually be \lq drawn\rq\  on the plane, 
so we may regard it as a figure on the plane.
We call this figure a {\em diagram of $L$}. 
In general, a {\em link diagram} is a union of closed curves on the plane with finitely many crossings, 
where each crossing  has an \lq over/under information\rq.
Likewise, we can define a {\em string link diagram} as a union of curves on the rectangle 
$D^1 \times [0,1]$ such that the $i$th component runs from $p_i\times\{0\}$ to $p_i\times\{1\}$. 
A crossing is locally the intersection of two orthogonal line segments, and 
the point with over (resp. under) information is called 
{\em over crossing} (resp. {\em under crossing}), where formally we assume that 
over and under crossings are distinct points on distinct line segments, 
see Figure~\ref{crossing}.

\begin{figure}[!h]
\begin{center}
\includegraphics[width=.85\linewidth]{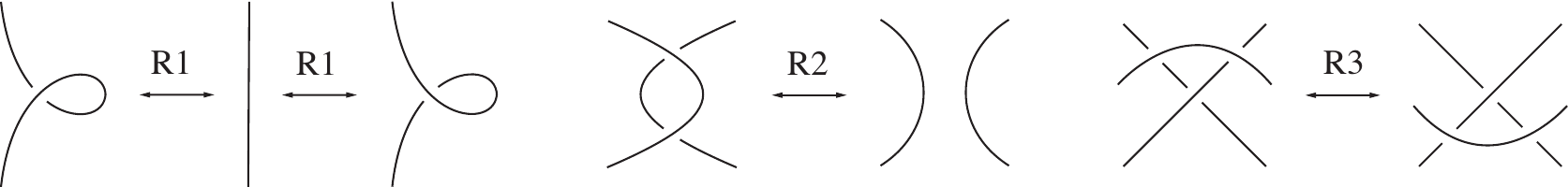}
   \vspace{.5em}
\caption{Reidemeister moves}\label{R123}
\end{center}
\end{figure}

Two (string) link diagrams are {\em equivalent} if one is deformed into the other 
by a combination of continuous deformations (fixing the boundary) and the three local moves 
R1, R2, R3 in Figure~\ref{R123}, called {\em Reidemeister moves}.

The following theorem shows that the classification of links is 
nothing more than the classification of link diagrams.

\medskip
\begin{thm} (Reidemeister Theorem)
Two links are equivalent if and only if 
their diagrams are equivalent. 
\end{thm}

\medskip
\subsection{Virtual link diagrams and welded links {\rm (\cite{Kau},\cite{GPV},\cite{FRR})}}

While all crossings of link diagrams have over/under informations, 
we extend the concept of link diagrams by considering diagrams that formally allow crossings  
without over/under informations as on the left side of Figure~\ref{crossing}. 
A diagram containing crossings without over/under information does not correspond to any 
link in 3-space. It is an \lq imaginary\rq\ link, so it is called a {\em virtual link diagram}.
A crossing without over/under information is called a {\em virtual crossing}. 
On the other hand, a crossing with an over/under information is called a {\em classical crossing}. 
For convenience, we also allow the case where there are no virtual crossings in a virtual link diagram.
From now on, when we emphasize that link diagrams contain only classical crossings, we call 
them  {\em classical link diagrams}.

\begin{figure}[htbp]
\begin{center}
\includegraphics[width=.7\linewidth]{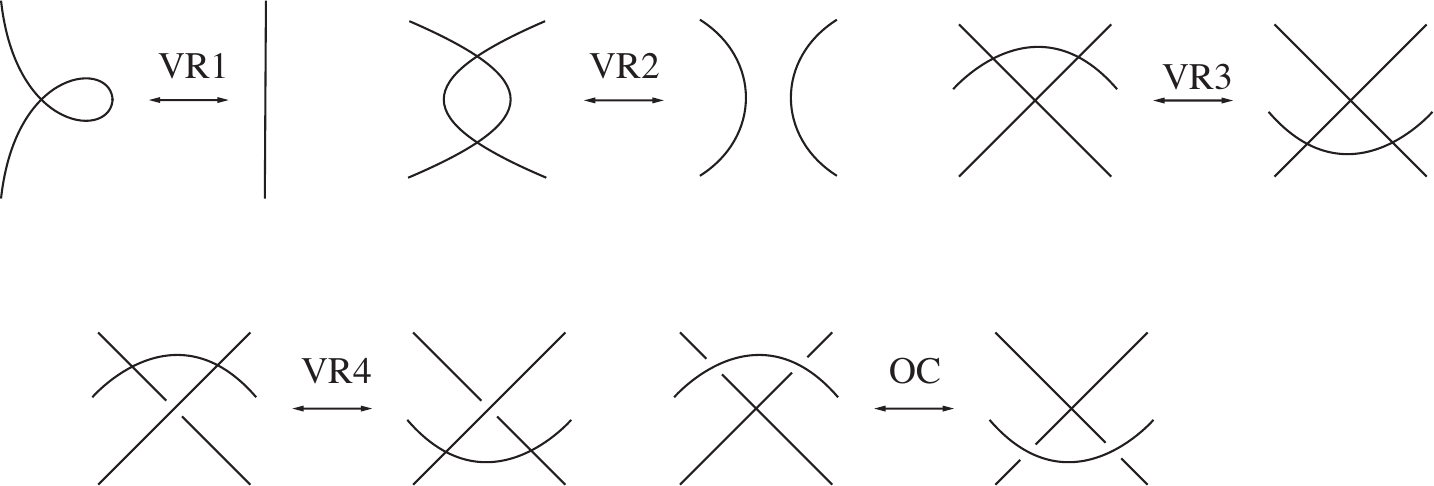}
   \vspace{.5em}
\caption{Virtual Reidemeister moves, OC move}\label{WR}
\end{center}
\end{figure}

For virtual link diagrams, we define five local moves as follows. 
The first four moves in Figure~\ref{WR} are called 
{\em virtual Reidemeister moves} (VR1,VR2,VR3,VR4), and 
the last move is called {\em OC move} ({\bf O}vercrossings {\bf C}ommute).

The residue classes of the set of virtual link diagrams modulo 
continuous deformations, Reidemeister moves, 
virtual Reidemeister moves and OC moves is called 
{\em welded links}.\footnote{The residue classes modulo continuous deformations, Reidemeister moves 
and virtual Reidemeister moves
are called {\em virtual links}. Here we never treat virtual links.}
We define {\em welded string links} in a similar way. 

The Reidemeister moves, the virtual Reidemeister moves, and the OC move are 
called the {\em welded moves}.
As an important property of welded links, the following theorem is known.

\medskip
\begin{thm} (\cite{T},\cite{GPV})
Two classical link diagrams are equivalent as classical links if and only if  
they are equivalent as welded links. More precisely, two classical link diagrams are 
equivalent up to continuous deformations and Reidemeister moves if and only if 
they are equivalent up to continuous deformations and welded moves.
\end{thm}

\medskip
This theorem implies that there is a natural injection 
from the equivalence classes of classical links to 
welded links, i.e., we can say 
\[\text{\lq classical links can be embedded in 
welded links.\rq}\footnote{The author does not know if string links can be embedded in 
welded string links.}\]
In this sense, as we mentioned in introduction, 
the relationship between \lq classical links\rq\ and \lq welded links\rq\  is similar to 
the relationship between \lq real numbers' and \lq complex numbers\rq. 

\medskip
\subsection{Based welded links}

\begin{figure}[!h]
\begin{center}
\includegraphics[width=.3\linewidth]{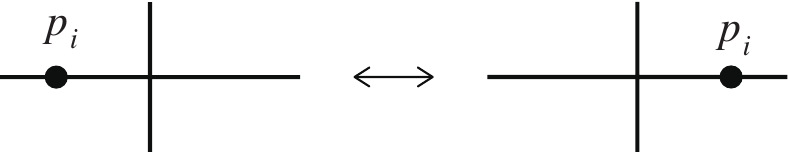}
   \vspace{.5em}
\caption{A base point passing through a virtual crossing}\label{BC}
\end{center}
\end{figure}

For an $n$-component virtual link diagram $L$, we choose $n$ points $p_1,...,p_n$ on 
$L$ such that each $p_i$ is on the $i$th component of $L$ and disjoint from the crossings of $L$. 
The pair $(L,\mathbf{p})$ of $L$ and $\mathbf{p}=(p_1,...,p_n)$ is called a {\em based virtual link diagram}.
Two based virtual link diagrams are {\em equivalent}  
 if one is deformed into the other 
by a combination of continuous deformations, moves as illustrated in Figure~\ref{BC}, 
and welded moves that do not contain base points.  
The residue classes of the set of based virtual link diagrams modulo 
this equivalence relation is called {\em based welded links}. 
Based welded links are \lq between\rq\  welded string links and welded links. 
In fact we have the following sequence of surjections:
\[\{\text{welded string links}\}\twoheadrightarrow
\{\text{based welded links}\}\twoheadrightarrow
\{\text{welded links}\}.\]
One advantage of considering based welded links is that peripheral systems  
can be uniquely defined (as explained in the next section).
This means that, like string links, based welded links
are very \lq suitable\rq\  objects for defining Milnor invariants.
In the next section, we define Milnor invariants for based welded links, and 
by using this definition, we define Milnor invariants for welded string links and 
for welded links.

\medskip
\subsection{Peripheral systems of based welded links}

A {\em diagram} of a (based/string) welded link means a representative of the welded link, i.e., a virtual link diagram. 
However, since it is tedious to distinguish them, the welded link (i.e. the residue class) and its diagram 
(i.e. a representative) are often identified.
From now on we simply call based welded link diagrams {\em based diagrams}.

\begin{figure}[!h]
\begin{center}
\begin{overpic}[width=9cm]{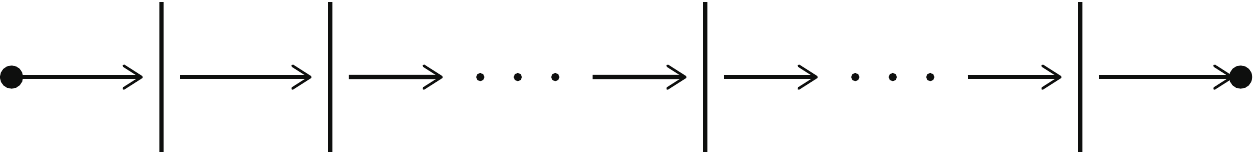}
\put(0,5){$p_{i}$}
\put(13,22){$a_{i0}$}
      \put(45,22){$a_{i1}$}
      \put(75,22){$a_{i2}$}
      \put(110,22){$a_{i(j-1)}$}
      \put(150,22){$a_{ij}$}
      \put(185,22){$a_{i(r_{i}-1)}$}
      \put(230,22){$a_{ir_{i}}$}
      \put(30,-7){$u_{i1}$}
      \put(65,-7){$u_{i2}$}
      \put(141,-7){$u_{ij}$}
      \put(218,-7){$u_{ir_{i}}$}
\put(250,5){$p_{i}$}
    \end{overpic}
  \end{center}
   \vspace{.5em}
  \caption{Arcs of $K_i$}
  \label{schematic}
\end{figure}
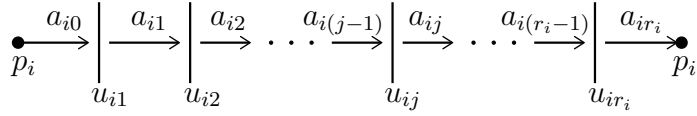

For a based diagram $(L, \mathbf{p})$, $L$ is divided into several segments by its base points and under-crossings. 
These segments are called the {\em arcs} of $L$.
Let $a_{i0}$ be the outgoing arc from the base point $p_i$, and 
let $a_{i1},...,a_{i{r_i}}$ be the other arcs of the $i$th component of $L$ 
that appear after $a_{i0}$ in order, when traveling around the $i$th component from $p_i$ along the orientation, 
where $r_i+1$ is the number of arcs of the $i$th component 
($i=1,...,n)$, see Figure~\ref{schematic}. Note that $a_{ir_i}$ is the ingoing arc to $p_i$

The {\em group $G(L,\mathbf{p})$ of $(L,\mathbf{p})$} is the quotient group of the 
free group $\widetilde{F}$ with generating set $\{a_{ij}\}_{i,j}$ modulo 
the following relations
\[R_{ij}=a_{i(j-1)}u_{ij}^{\e(ij)}a_{ij}^{-1}u_{ij}^{-\e(ij)}~~(1\leq i\leq n,~1\leq j\leq r_i),\]
where $u_{ij}$ denote the arc containing the over crossing between the arcs $a_{i(j-1)}$ and $a_{ij}$ 
as illustrated in Figure~\ref{schematic}, 
and $\varepsilon(ij)\in\{-1,1\}$ is the {\em sign of the crossing} involving   $a_{i(j-1)},a_{ij}$ and $u_{ij}$;  
in the figure, if the orientation of $u_{ij}$ is from up to down, then $\varepsilon(ij)=+1$, 
otherwise $\varepsilon(ij)=-1$.

\medskip
\begin{remark}\label{group}
We note that $G(L,\mathbf{p})$ has a group presentation 
$\langle \{a_{ij}\}_{i,j}~|~ \{R_{ij}\}_{i,j} \rangle$. 
Let $G(L)$ be the group obtained from $\langle \{a_{ij}\}_{i,j}~|~ \{R_{ij}\}_{i,j} \rangle$ 
by adding the relations $a_{i0}a_{ir_i}^{-1}~(i=1,...,n)$. 
We call $G(L)$ the {\em group of  $L$}. It is known that if $L$ is a classical link, then 
$G(L)$ is isomorphic to the fundamental group of the complement of $L$. 
\end{remark}

\medskip
For each $i$, the two elements $a_{i0}$ and 
\[\lambda_i=a_{i0}^ {-w_i}u_{i1}^{\varepsilon(i1)}u_{i2}^{\varepsilon(i2)}\cdots u_{ir_i}^{\e(ir_i)}
~(\text{where $w_i$ is the sum of signs $\e(il)$ for all $u_{il}\subset K_i$} )\]
of $G(L,\mathbf{p})$ are called the {\em $i$th meridian} and the {\em $i$th longitude} respectively. 
As in the case of classical links, the pair 
$(G(L,\mathbf{p}),\{a_{i0},\lambda_i\}_i)$ is called the {\em peripheral system} of $(L,\mathbf{p})$. 
We note that the peripheral system of $G(L,\mathbf{p})$ is uniquely determined from 
the diagram $(L,\mathbf{p})$.


\medskip
\subsection{Peripheral systems of diagrams and Milnor invariants}\label{peri-diagram}

We consider the natural projection 
\[\rho_q:G(L,\mathbf{p}) \longrightarrow N_q(L,\mathbf{p})=G(L,\mathbf{p})/\Gamma_qG(L,\mathbf{p}).\]
We remark that $N_q(L,\mathbf{p})$ is a nilpotent group  
generated by $\rho_q(a_{i0})=\alpha_i~(i=1,...,n)$, and 
it is called the {\em $q$th nilpotent group} of $(L,\mathbf{p})$ (or the {\em $q$th nilpotent quotient} of $G(L,\mathbf{p})$).  
Hence each $\rho_q (\lambda_k)$ can be written as a word of 
$\alpha_1,...,\alpha_n$. 
Moreover, it is known that 
for the free group $F$ with generating set $\{\alpha_i~|~(i=1,...,n)\}$, we have
\[N_q(L,\mathbf{p})\cong F/\Gamma_qF, \]
see Theorem~\ref{CM}.
As in the case of string links (see the last paragraph in Section 1.2), 
by taking the Magnus expansion $E(\rho_q(\lambda_k))$, we have 
integers $\mu_{(L,\mathbf{p})}(i_1,...,i_{s} k)~(s<q)$, which are invariants for 
the based diagram $(L,\mathbf{p})$.
Since the condition $s<q$ is not essential, we have invariants for all sequences.  
We call  them  {\em Milnor $\mu$-invariants for based welded links}.  
 
For a welded string link $S$, by identifying both ends $p_i\times\{0\},p_i\times\{1\}$ of $S$ 
to a point $p_i~(i=1,...,n)$, we obtain a based diagram $(L,\mathbf{p})$.
Then we define $\mu_S(I):=\mu_{(L,\mathbf{p})}(I)$ for all sequences $I$, and 
call them  {\em Milnor $\mu$-invariants for the welded string link $S$}.
For a welded link $L$, by choosing a set $\mathbf{p}$ of base points  on $L$, we have 
a based diagram $(L,\mathbf{p})$. 
And we define {\em Milnor $\omu$-invariants $\omu_L(I)$ for the welded link $L$} for all sequences $I$ 
as the residue classes of $\mu_{(L,\mathbf{p})}(I)$ modulo $\Delta_{(L,\mathbf{p})}(I)$, where 
\[\Delta_{(L,\mathbf{p})}(I):=\gcd
\left\{\mu_{(L,\mathbf{p})}(J)~\vline~
\begin{array}{l}
J\textrm{ is a sequence obtained from $I$ by deleting at least }\\
\textrm{one index and permuting the resulting sequence cyclicly}
\end{array}
\right\}.\]
As the names \lq invariants\rq\  suggest, Milnor invariants for welded string links and  
welded links are invariants under their equivalence relations respectively.


\medskip
\subsection{Automorphisms of nilpotent groups of diagrams and Milnor invariants}\label{aut-nil-diagram}

For each $k(=1,...,n)$, since  
$\rho_q(a_{kr_k})=(\rho_q (\lambda_k))^{-1} \alpha_k\rho_q (\lambda_k)$, we have 
the automorphism
\[\varphi_q=\varphi_q(L,\mathbf{p}):
F/\Gamma_q F\longrightarrow F/\Gamma_q F,~\alpha_k\longmapsto 
\rho_q(a_{kr_k}).\]
By applying the five lemma to the natural exact sequence 
\[1\rightarrow\Gamma_{q-1} F/\Gamma_q F\hookrightarrow
F/\Gamma_q F\twoheadrightarrow F/\Gamma_{q-1} F\rightarrow 1, \]
we can see that 
$\varphi_q(L,\mathbf{p})$ is an isomorphism. 
As in the case of classical string links, $\varphi_q(L,\mathbf{p})$ can be regarded as  
the Milnor invariants of $(L,\mathbf{p})$ of length at most $q-1$.
In fact, by using Proposition~\ref{Magnus}, we have the following proposition.  

\medskip
\begin{prop}\label{meilhan-yasu-lemma}(cf. \cite[Remark~3.6]{MY2})
Let $q$ be a positive integer, and let $(L,\mathbf{p})$ and $(L',\mathbf{p}')$ 
be based diagrams. 
For any sequence $I$ of length at most $q-1$, 
$\mu_{(L,\mathbf{p})}(I)=\mu_{(L',\mathbf{p}')}(I)$ if and only if 
$\varphi_{q}(L,\mathbf{p})=\varphi_{q}(L',\mathbf{p}')\in 
\mathrm{Aut}_{\mathrm{c}}(F/\Gamma_{q} F)$, 
where 
$\mathrm{Aut}_{\mathrm{c}}(F/\Gamma_{q}F)$ is 
the set of automorphisms of $F/\Gamma_{q} F$
that act by conjugation on each generator. 
\end{prop}

\medskip
\subsection{Colorings of diagrams and Milnor invariants}\label{color-diagram}

A map $f$ from the set of arcs of a diagram $(L,\mathbf{p})$ to 
a group $X$  is called an {\em $X$-coloring} of $(L,\mathbf{p})$ 
if $f$ satisfies 
\[f(a_{ij})=f(u_j)^{-\e(j)}f(a_{i(j-1)})f(u_j)^{\e(j)},\]
for each classical crossing involving $a_{i(j-1)},a_{ij}$ and $u_j$ (see Figure~\ref{schematic}). 
Let $(L',\mathbf{p}')$ be a diagram obtained from $(L,\mathbf{p})$ by a single welded move. 
Then it is easily seen that there is an $X$-coloring of $(L',\mathbf{p}')$ which coincides 
with $f$ for the arcs not contained in  the disk where the welded move is applied.  

Let $\widetilde{N}$ be the normal closure of $\widetilde{F}$ that contains  
$\{R_{ij}\}_{i,j}$, and let $A(L,\mathbf{p})$ be the set of arcs $\{a_{ij}\}_{i,j}$ of $(L,\mathbf{p})$.
Then by compositing the natural projection 
\[A(L,\mathbf{p})\subset \widetilde{F}\twoheadrightarrow\widetilde{F}/\widetilde{N}
= G(L,\mathbf{p})\]
 and $\rho_q$, we have the map 
\[\phi_q:A(L,\mathbf{p}) \longrightarrow F/\Gamma_q F(\cong N_q(L,\mathbf{p})).\]
Since $\phi_q$ is an $(F/\Gamma_q F)$-coloring of $(L,\mathbf{p})$ and uniquely determined 
by $(L,\mathbf{p})$, $\phi_q(a_{kr_k})$ is an invariant for $(L,\mathbf{p})$. 
Moreover since 
$\varphi_q(\alpha_k)=\phi_q(a_{kr_k})$,  $\phi_q$ determines 
the automorphism $\varphi_q$ of the nilpotent group of $(L,\mathbf{p})$.
By Proposition~\ref{meilhan-yasu-lemma}, these $(F/\Gamma_q F)$-colorings can also be regarded as
Milnor invariants of length at most $q-1$. 

\medskip
\subsection{Chen-Milnor map and an algorithm for computing Milnor invariants}

Here we show how to calculate Milnor invariants of based diagrams.
This is simply a rewriting of the method given by Milnor \cite{M2} for classical links into one for based diagrams. 

Let $(L,\mathbf{p})$ be a based diagram. For the $i$th component $K_i$ illustrated in Figure~\ref{schematic}, 
we put  
\[v_{ij}=u_{i1}^{\varepsilon(i1)}\cdots u_{ij}^{\varepsilon(ij)}.\]
(Note that $\lambda_i=a_{i0}^ {-w_i}v_{ir_i}$.)
Let $F$ be the free group generated by $a_{i0}=\alpha_i~(i=1,...,n)$. 
Then, for a positive integer $q$, we define inductively a homomorphism  
\[
\eta_{q}=\eta_{q}(L,\mathbf{p}):\widetilde{F}\longrightarrow F(\subset \widetilde{F})\] 
as follows:\footnote{While $\eta_q$ is essentially the same as defined in \cite{M2}, 
in \cite{M2} $a_{i0}\cup a_{ir_i}$ is treated as a single arc.}
\[
(\mathrm{i})~\eta_{1}(a_{ij})=\alpha_{i};~~~
(\mathrm{ii})~\eta_{q+1}(a_{i0})=\alpha_{i},~\eta_{q+1}(a_{ij})=\eta_{q}(v_{ij}^{-1})\alpha_{i}\eta_{q}(v_{ij}) \hspace{1em} (1\leq j\leq r_i). \]
The map $\eta_q$ is called {\em Chen-Milnor map} \cite{Chen},\cite{M2}. 
The following lemma can be easily shown by induction on $q$.

\medskip
\begin{lem}({\cite{Chen},\cite{M2}})\label{Chen}
\begin{enumerate}
\item For each element $a_{ij}\in\widetilde{F}$,\\
~~~~~{\rm(i)}~$\eta_{q}(a_{ij})\equiv a_{ij} \pmod{\Gamma_q\widetilde{F}\cdot\widetilde{N}}$, and 
~~~{\rm(ii)}~$\eta_{q}(a_{ij})\equiv\eta_{q+1}(a_{ij}) \pmod{\Gamma_{q}F}$. 
\item For each relation $R_{ij}$, $\eta_{q}(R_{ij})\equiv 1\pmod{\Gamma_{q}F}$. 
\end{enumerate}
\end{lem}

\medskip
Since $N_q(L,\mathbf{p})=(\widetilde{F}/\widetilde{N})/\Gamma_q(\widetilde{F}/\widetilde{N})
\cong \widetilde{F}/\Gamma_q\widetilde{F}\cdot\widetilde{N}$, by the lemma above, we have 
the following theorem. 

\medskip
\begin{thm}(J. Milnor \cite{M2})\label{CM}~
\begin{enumerate}
\item $\eta_q$ induces an isomorphism
\[N_q(L,\mathbf{p})\cong\langle\alpha_1,\ldots,\alpha_n~|~ \Gamma_qF\rangle.\]
\item $\rho_q(\lambda_k)=\eta_q(\lambda_k){\Gamma_qF}
(\in F/\Gamma_{q} F)$, and hence 
$\eta_q(\lambda_k)$ is a representative of $\rho_q(\lambda_k)$.
\end{enumerate}
\end{thm}

\medskip
\begin{remark}
Since $\eta_q(a_{i0}a_{ir_i}^{-1})=[\alpha_i,\alpha_i^{w_i}\eta_q(\lambda_i)]=[\alpha_i,\eta_q(\lambda_i)]$ for each $i(=1,...,n)$, 
by combining Remark~\ref{group} and 
Theorem~\ref{CM}, we see that $N_q(L)=G(L)/\Gamma_qG(L)$ has the following presentation
\[\langle\alpha_1,\ldots,\alpha_n~|~
[\alpha_1,\eta_q(\lambda_1)],\ldots ,
[\alpha_n,\eta_q(\lambda_n)], \Gamma_qF\rangle.\]
\end{remark}

\medskip
By Proposition~\ref{Magnus}~(2) and Theorem~\ref{CM}, we have the following corollary. 

\medskip
\begin{cor}
The Milnor invariant $\mu_{(L,\mathbf{p})}(i_1,...,i_sk)$ for $(L,\mathbf{p})$ is 
equal to the 
coefficient of $X_1\cdots X_s$ in the Magnus expansion of 
$\eta_q(\lambda_k)$, where 
$q$ is an integer with $q>s$.
\end{cor}

\medskip
\begin{remark}
As explained above, Milnor invariants are \lq theoretically\rq\  easy to compute.
However, when we actually do the calculations, we notice that the length of the 
word $\eta_q(\lambda_k)$ grows exponentially with $q$.
Furthermore,  performing the Magnus expansion on it, the amount of 
calculations becomes enormous.
Therefore, by this naive algorithm, it is impossible to calculate it even using a computer.
Takabatake, Kuboyama and Sakamoto noticed that long words 
that appear in the process of calculating Milnor invariants contain a lot of repetitions, 
and by \lq folding\rq\ the words by these repetitions, they succeeded to produce a program 
that can calculate up to about $q=16$ (depending on the complexity of the links) \cite{TKS}. 
\end{remark}

\medskip
\begin{remark}
Dye and Kauffman's 2010 paper is often cited as the first paper that tried to extend Milnor invariants to virtual diagrams. However, 
their paper contains obvious and fatal errors (\cite[Remark~4.6]{Kotorii}, \cite[Remark 6.8]{MWY}).
Therefore it should be avoided. 
In fact, the first successful extension is by Kravchenko and Polyak in~\cite{KP}. 
They extended Milnor link-homotopy $\mu$-invariants to virtual string links. 
Kotorii then extended Milnor link-homotopy $\omu$-invariants to virtual links in~\cite{Kotorii}.  
Both extensions are invariants of welded diagrams, but restricted to the case of link homotopy invariants. 
In the general case,  
Audoux, Bellingeri, Meilhan and Wagner extended Milnor $\mu$-invariants to welded string links in~\cite{ABMW}, 
and 
Chrisman defined Milnor $\omu$-invariants for welded links in~\cite{C}. 
In \cite{MWY}, Miyazawa, Wada and the author gave Milnor invariants for (based) welded (string) links 
using a way completely independent of \cite{ABMW} and \cite{C}.
Although these extended Milnor invariants are defined in different ways,
they are the same invariant, i.e., for welded (string) links, their values coincide.
\end{remark}

\section{Characterization of Milnor invariants}

In this chapter we discuss the characterization of based diagrams that have the same Milnor invariants.
Milnor invariants can be characterized by two equivalence relations, {\em $W_k$-concordance} and 
{\em self $W_k$-concordance}, on diagrams. 

\medskip
\subsection{$W_k$-concordance and Milnor invariants}
The $W_k$-concordance is an equivalence relation that combines two equivalence relations, 
{\em $W_k$-equivalence} and {\em welded-concordance}, which are explained in the next section. 
The $W_k$-equivalence can be seen as a welded version of $C_k$-equivalence, 
an equivalence relation on classical links defined by Habiro~\cite{H}.

B. Colombari \cite{Col} showed the following theorem.

\medskip
\begin{thm}\label{Boris}(B. Colombari \cite{Col})
Let $k$ be a positive integer. 
Two based diagrams
$(L,\mathbf{p})$ and $(L',\mathbf{p}')$ are $W_k$-concordant if and only if 
$\mu_{(L,\mathbf{p})}(I)=
\mu_{(L',\mathbf{p}')}(I)$ for any sequence $I$ of length at most $k$. 
\end{thm}

\medskip
\begin{remark}
In the paper \cite{Col}, the statement of Theorem~\ref{Boris} is made for welded string links 
rather than based diagrams. These statements are essentially the same.
We can define a {\em product} on the set $\mathrm{wSL}$ of welded $n$-string links. 
It is shown that the set of $W_k$-concordance classes $\mathrm{wSL}/{(W_k+\mathrm{c})}$  
forms a group under the product. 
Moreover, by using Proposition~\ref{meilhan-yasu-lemma} (for the welded string link version) and 
Theorem~\ref{Boris}, we can show the following isomorphism
\[\mathrm{wSL}/{(W_k+\mathrm{c})} \cong \mathrm{Aut}_{\mathrm{c}}(F/\Gamma_{k+1} F).\]
\end{remark}

\medskip
\begin{remark}
Theorem~\ref{Boris} also holds when restricted to classical string links.
That is, it completely characterizes Milnor invariants of classical string links.
There are many papers that try to characterize Milnor invariants.
In particular, the {\em $C_k$-concordance} \cite{MY0}  and {\em Whitney tower concordance} \cite{CST} are 
closely related to Milnor invariants, 
but they cannot give a complete characterization.
These studies are done in the \lq real world\rq\  of classical string links and $C_k$-concordance.
On the other hand, Colombari gave the result in the \lq real world\rq\  by expanding his study to 
the \lq imaginary world\rq\  (i.e., welded string links and $W_k$-concordance).
This is a successful example of how extending the world of classical objects to include welded objects 
allows us to solve problems from a new perspective.
\end{remark}

 \medskip
\subsection{$W_k$-equivalence}

The content of this section is due to \cite{MY1}.
Our main purpose here is to introduce the definition of $W_k$-equivalence.
Before explaining $W_k$-equivalence, some preparation is necessary.
From now on, in order to deal with diagrams of (based) welded links and welded string links simultaneously, 
for convenience, we assume that they are contained in the disk.

\medskip 
\subsubsection{$W$-tree}
Let $L$ be a diagram.  
A tree $T$ that is \lq drawn\rq\ ({\em immersed}) on the disk is called 
a $W$-tree for  $L$ if it satisfies the following conditions:
\begin{enumerate}
\item  The valency of each vertex of $T$ is either $1$ or $3$, i.e., $T$ is a {\em uni-trivalent tree}.
\item The trivalent vertices are disjoint from $L$, and the univalent vertices are contained in $L$ 
and disjoint from the crossings (and base points) of $L$.
\item All edges are oriented so that 
each trivalent vertex has two ingoing and one outgoing edge. 
\item We allow only virtual crossings between edges of $T$, and between $L$ and edges of $T$.
\item Each edge of $T$ is assigned a number (possibly zero) of decorations, called
{\em twists}, which are disjoint from all vertices and crossings.
(Although base points and twists are both presented by black dots, 
twists are on edges of $T$ and base points are on the diagram $L$.) 
\end{enumerate}

A univalent vertex of $T$ with outgoing edge is called a {\em tail}, and 
a univalent vertex with an ingoing edge is called a {\em head}. 
From condition (3) above, we can see that a $W$-tree has a unique head, and that 
specifying the head determines the tails and the orientation of the entire $W$-tree.
The tails and head are called {\em ends}.

The {\em degree of $T$} is the number of the tails of $T$, and a
{\em $W_{k}$-tree} is a $W$-tree with degree $k$. 
In particular, a ${W}_{1}$-tree is called a {\em $W$-arrow}.

\begin{figure}[!h]
  \begin{center}
\includegraphics[width=.3\linewidth]{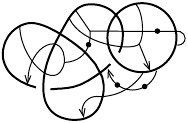}
   \vspace{.5em}
    \caption{A union of a $W_3$-tree, a $W_2$-tree and two $W$-arrows}\label{extree}
  \end{center}
\end{figure}

For a union of $W$-trees for $L$, vertices are assumed to be pairwise disjoint, and
all crossings among edges are assumed to be virtual, see Figure~\ref{extree} for an example.
From now on, when drawing diagrams, we will distinguish them from $W$-trees 
by drawing them with slightly thicker lines as in Figure \ref{extree}.

\medskip
$W$-trees act as \lq guides\rq\ for transforming diagrams locally.
{\em Surgery along a $W$-tree} is a local move on a diagram 
using the $W$-tree as a guide, as described below.
Before describing surgery along a $W$-tree, we need to explain {\em surgery along $W$-arrow}.

\subsubsection{Surgery along $W$-arrows}
Let $A$ be a union of $W$-arrows for $L$. 
{\em Surgery along $A$} yields a new diagram, denoted by $L_A$, 
which is defined as follows.

\begin{figure}[!h]
  \begin{center}
    \begin{overpic}[width=9cm]{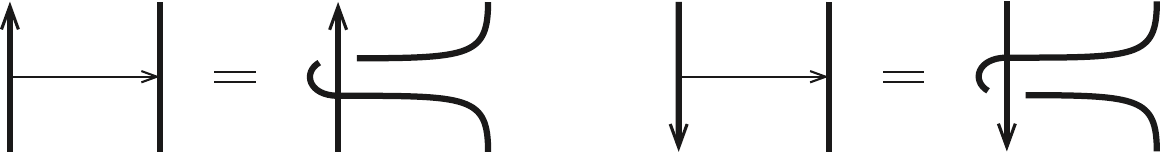}
      \put(5,-12){$L\cup A$} 
      \put(85,-12){$L_{A}$} 
      \put(152.5,-12){$L\cup A$}
      \put(232,-12){$L_{A}$} 
    \end{overpic}
  \end{center}
   \vspace{.5em}
  \caption{Surgery along $W$-arrow~(i)}
  \label{surgery}
\end{figure}

\begin{figure}[!h]
  \begin{center}
    \begin{overpic}[width=12cm]{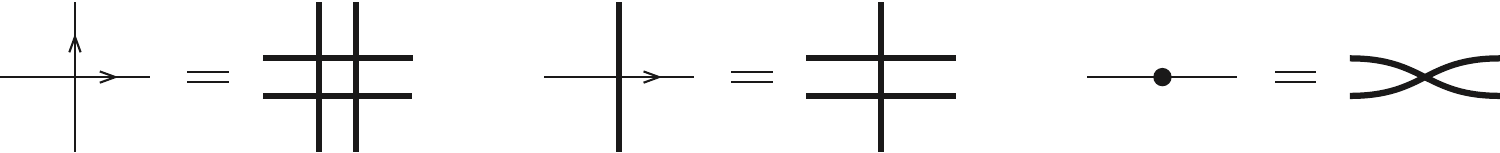}
      \put(-5,22){$A$}
      \put(14.5,-12){$A$} 
      \put(77,-12){$L_{A}$} 
      \put(122,22){$A$}
      \put(140,-12){$L$}
      \put(198,-12){$L_{A}$} 
      \put(259,-12){$A$}
      \put(317,-12){$L_{A}$} 
    \end{overpic}
  \end{center}
   \vspace{.5em}
  \caption{Surgery along $W$-arrow~(ii)}
  \label{surgery2}
\end{figure}
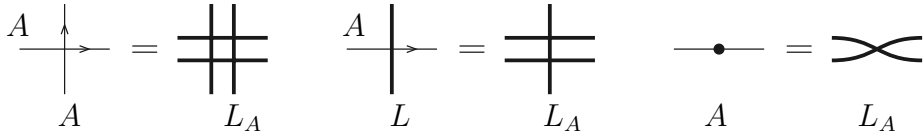

\begin{enumerate}
\item[(i)] If a $W$-arrow in $A$ does not contain any crossing or twist, 
we perform the operation shown in Figure~\ref{surgery}.
Note that this operation depends only on the orientation of $L$ near the tail of the $W$-arrow.
\item[(ii)] If a $W$-arrow in $A$ contains some crossings and/or twists, 
in addition to the operations in Figure~\ref{surgery} on the ends, 
we also perform the operations on the edge in Figure~\ref{surgery2}.
\end{enumerate}

\begin{figure}[!h]
  \begin{center}
    \begin{overpic}[width=6cm]{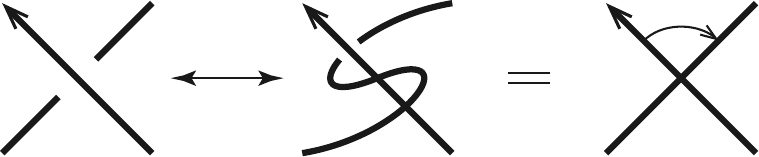}
      \put(42,22){VR2}
    \end{overpic}
  \end{center}
  \caption{A classical crossing can be replaced by 
  a virtual crossing and a $W$-arrow}
  \label{Aprst}
\end{figure}

An {\em arrow presentation} for a diagram $L$ is a pair $(V,A)$ of a
diagram $V$ without classical crossings and a union of $W$-arrows $A$ for $V$, 
such that $V_A$ is equivalent to the diagram $L$.
Any diagram admits an arrow presentation since  
all classical crossings can be replaced with virtual crossings and $W$-arrows as shown in Figure~\ref{Aprst}. 

Two arrow presentations $(V,A)$ and $(V',A')$ are {\em equivalent} if
$V_{A}$ and $V'_{A'}$ are equivalent. 
An arrow presentation $(V,A)$ is sometimes treated as a union $V\cup A$ of a diagram $V$ and $W$-arrows $A$, 
and we do not distinguish them unless necessary.
One of the advantages of the arrow presentation is that $V\cup A$ contains no classical crossings. 

\medskip
\subsubsection{Arrow-moves}
\begin{figure}[!h]
  \begin{center}
    \begin{overpic}[width=11.5cm]{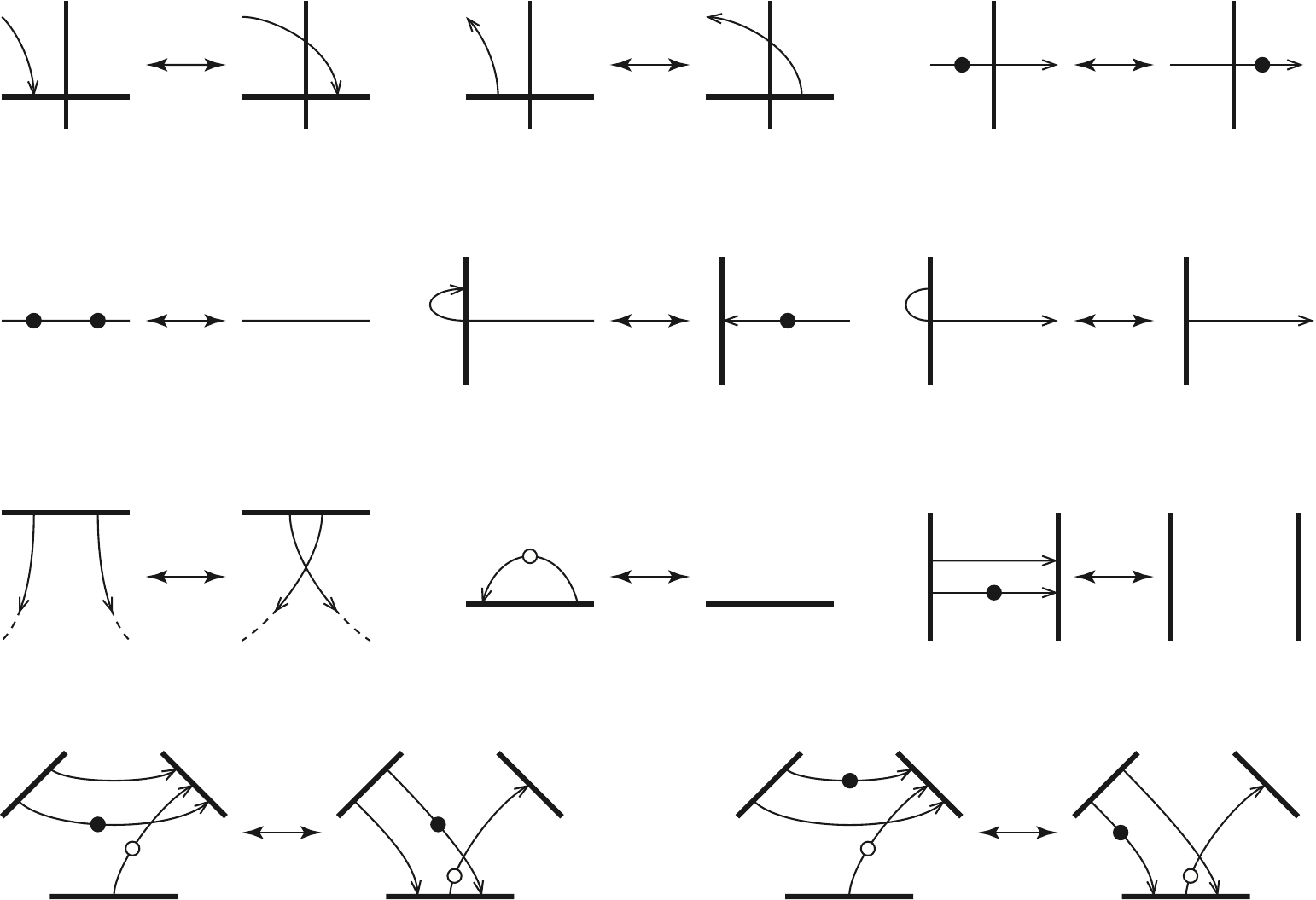}
      \put(37,212){AR1} 
      \put(152,212){AR2} 
      \put(267,212){AR3}
      \put(37,148.5){AR4} 
      \put(152,148.5){AR5}
      \put(267,148.5){AR6} 
      \put(37,84.5){AR7} 
      \put(152,84.5){AR8} 
      \put(267,84.5){AR9}
      \put(58,21){AR10} 
      \put(241.5,21){AR10}
    \end{overpic}
  \end{center}
  \caption{Arrow moves AR1--AR10}
  \label{Amoves}
\end{figure}
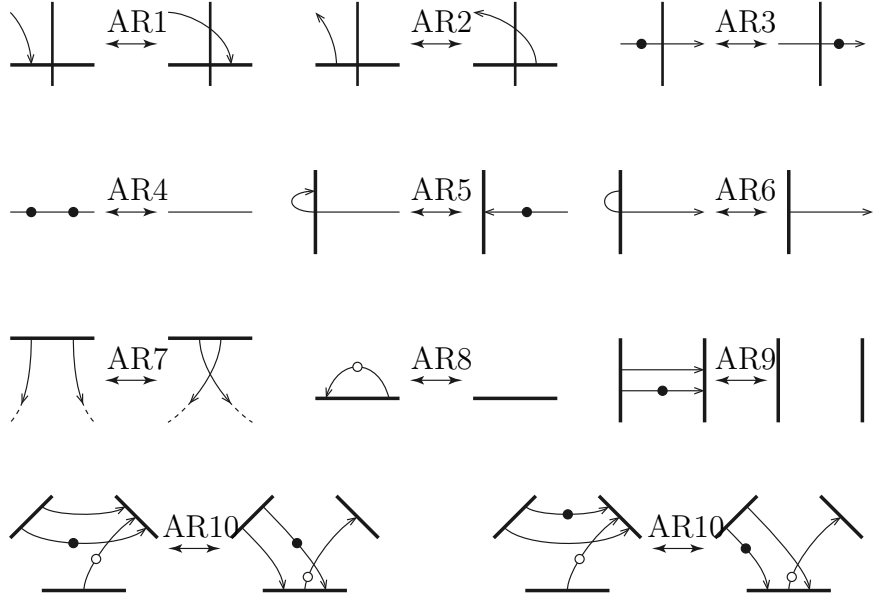

{\em Arrow moves} consist of the virtual Reidemeister moves VR1, VR2, VR3, 
where they may contain $W$-arrows, and 
the local moves AR1, ..., AR10 in Figure~\ref{Amoves}. 
We stress that arrow moves contain no classical crossings.
Here the vertical lines in AR1, AR2, AR3 are assumed to be included in diagrams or $W$-arrows, 
and white dots $\circ$ on $W$-arrows in AR8 and AR10 
mean that the $W$-arrows may or may not include twists.
Then we have the following. 

\medskip
\begin{thm}
Two arrow presentations are equivalent if and only if they are deformed into each other 
by a combination of continuous deformations and arrow moves, where 
in the case of based diagram, 
in addition to VR1, VR2, VR3, the move of Figure~\ref{BC} is also required, and 
the local moves AR1, ..., AR10 do not include base points.
\end{thm}

\medskip
\subsubsection{Surgery along $W$-trees}

\begin{figure}[!h]
  \begin{center}
    \begin{overpic}[width=11cm]{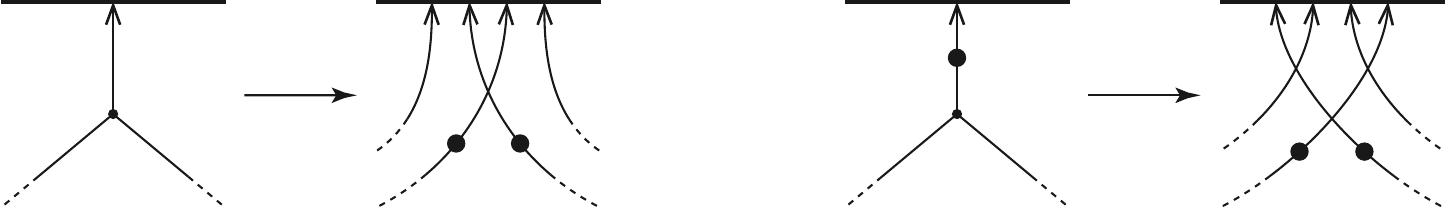}
      \put(58,29){(E)}
      \put(240,29){(E)}
    \end{overpic}
  \end{center}
  \caption{Expansion (E) of a tree}
  \label{expansion}
\end{figure}

\begin{figure}[!h]
  \begin{center}
    \begin{overpic}[width=12cm]{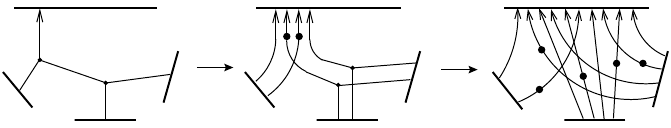}
      \put(103,37){(E)}
      \put(227,37){(E)}
    \end{overpic}
  \end{center}
  \caption{Expansions of a $W_3$-tree}
  \label{extension2}
\end{figure}

In the following, we define {\em surgery along W-trees}. 
For an integer $k\geq 2$, the expansion (E) of a $W_k$-tree is an operation 
that replaces a $W_k$-tree with four $W$-trees of degree less than $k$, 
as shown in Figure~\ref{expansion}.
For each figure of Figure~\ref{expansion}, the left (resp. right) dashed parts in the figure on the right of \lq$\stackrel{({\rm E})}{\longrightarrow}$\rq\ means the same figures as the left (resp. right) dashed part in the figure on the left 
of \lq$\stackrel{({\rm E})}{\longrightarrow}$\rq. 
The ends of the dashed parts are only the tails of the trees, and their arrangement is chosen arbitrarily.
In general, for a $W$-tree $T$, by repeating the expansion we obtain a union $E(T)$ of $W$-arrows from 
$T$ (Figure~\ref{extension2}).
While $E(T)$ depends on the arrangement of its tails, thanks to AR7 it is unique up to the arrow moves.
Then we define {\em surgery along a $W$-tree} $T$ 
as surgery along $E(T)$. In this case, we denote $L_{E(T)}$ by $L_T$.
Similarly, we define surgery along a union of $W$-trees. 

We say that $L_T$ is obtained from $L$ by a {\em $W_k$-move} if  $T$ is a $W_k$-tree.
In particular, we say that $L_T$ is obtained from $L$ by a {\em self $W_k$-move} 
if all ends of $T$ are contained in the same component of $L$.  
Using the move (I) in Figure~\ref{wtree-moves}, we can see that $L$ is obtained from $L_T$ by a 
(self) $W_k$-move, i.e., there is a (self) $W_k$-tree $T'$ such that $(L_T)_{T'}$ and $L$ 
are equivalent.

Two diagrams $L$ and $L'$ are {\em (self) $W_k$-equivalent} if there is a sequence of diagrams 
\[L=L_0,L_1,...,L_m=L'\]
such that, for each $i(=1,...,m)$, $L_i$ is equivalent to $L_{i-1}$, or 
$L_i$ is obtained from $L_{i-1}$ by a (self) $W_l$-move $(l\geq k)$.\footnote{It is shown that  
if $l>k$, then a (self) $W_l$-move is realized by (self) $W_k$-moves \cite{AMY1}.}

A {\em tree presentation} for a diagram $L$ is a pair $(V,T)$ of a
diagram $V$ without classical crossings and a union of $W$-trees $T$ for $V$, 
such that $V_T$ is equivalent to the diagram $L$.
Two tree presentations $(V,T)$ and $(V',T')$ are {\em equivalent} if
$V_{T}$ and $V'_{T'}$ are equivalent. 
A tree presentation $(V,T)$ is sometimes treated as a union $V\cup T$ of a diagram $V$ and $W$-trees $T$, 
and we do not distinguish them unless otherwise necessary.

\medskip
\subsubsection{$W$-tree moves}

A local move on a $W$-tree presentation is called a {\em $W$-tree move}
if it preserves the equivalence classes of the $W$-tree presentations.
Here, we introduce four of them given in \cite{MY1} that will be needed later.

\begin{figure}[!h]
  \begin{center}
    \begin{overpic}[width=7cm]{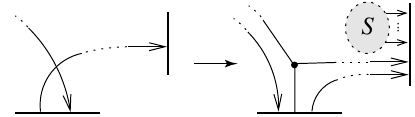}
      \put(80,-13){move (HTE)}
    \end{overpic}
  \end{center}
  \vspace{1em}
  \caption{$W$-tree moves}\label{wtree-moves}
 \label{HTE}
\end{figure}

The moves in Figure~\ref{wtree-moves} are called 
{\em move (I)} ({\bf I}nverse), {\em move (TE)} ({\bf T}ail {\bf E}xchange), 
{\em move (HE)} ({\bf H}ead {\bf E}xchange) and 
{\em move (HTE)} ({\bf H}ead {\bf T}ail {\bf E}xchange). 
 
In the move (TE), the two tails may be in the same $W$-tree. 
The right side $S$ of the move (HTE) in Figure~\ref{wtree-moves} represents a union of several $W$-trees, each of which 
has degree greater than the sum of the degrees of the two $W$-trees on the left side.
Moves (I) and (TE) are obtained by AR9 and by a combination of AR7 and the expansions of $W$-trees 
respectively. 
For the move (HE), by the expansion of the middle $W$-tree of the $W$-tree presentation on the right side 
and applying the move (I) and AR4, the $W$-tree presentation on the left side is obtained.
 To obtain the move (HTE), in addition to AR9, AR10, the expansions of $W$-trees and the move (HE), 
 a move called \lq Twist\rq\ (not introduced here) is required. 
 This requires more complex discussions than the other moves, 
 so we omit the details.
 
\begin{figure}[!h]
  \begin{center}
    \begin{overpic}[width=12.5cm]{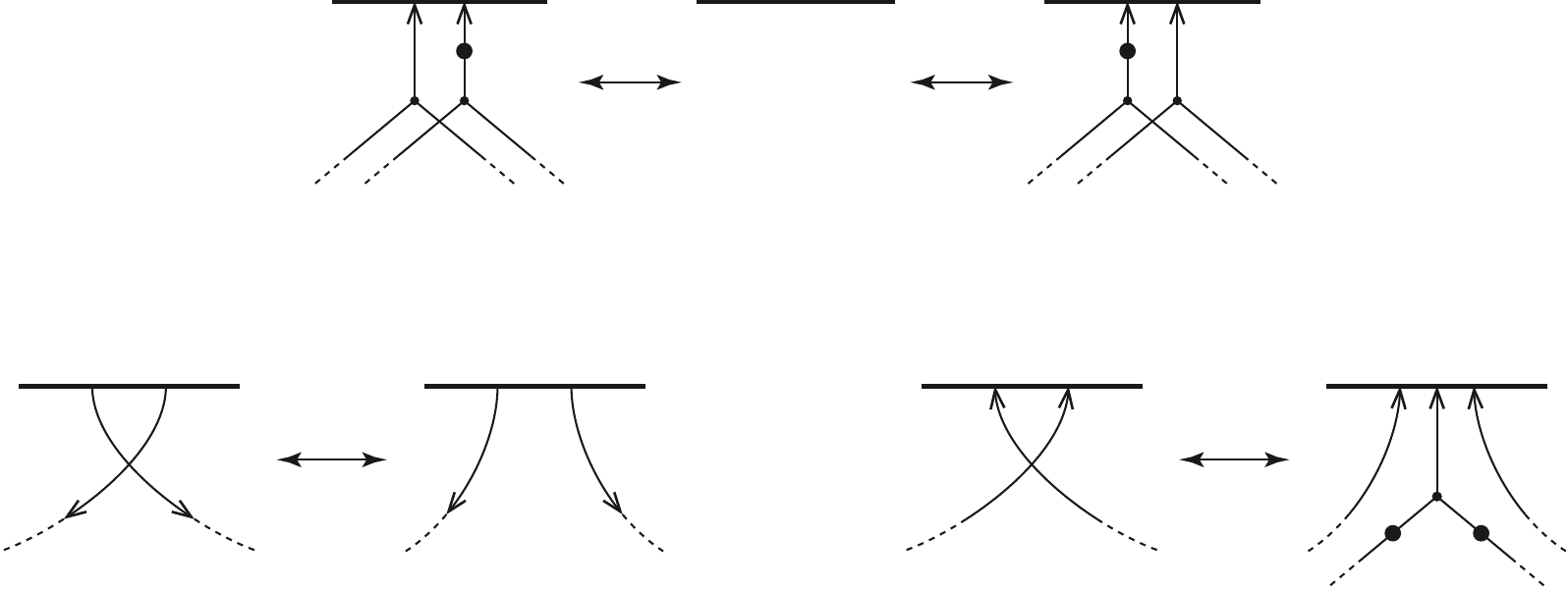}
      \put(165,76){move (I)}
      \put(58,-15){move (TE)}
      \put(264,-15){move (HE)}
    \end{overpic}
  \end{center}
\end{figure}

\medskip
\subsection{Welded-concordance}

\begin{figure}[!h]
\begin{center}
 \includegraphics[scale=0.9]{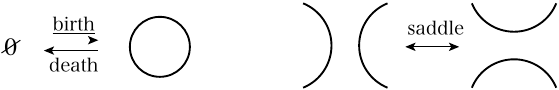}
   \vspace{1em}
 \caption{Birth/death and saddle moves}\label{fig:concmoves}
\end{center}
\end{figure}

Two $n$-component diagrams $L$ and $L'$ are \emph{welded-concordant} if one can be deformed 
into the other by a sequence of
welded equivalence and the {\em birth/death} and {\em saddle} moves of Figure~\ref{fig:concmoves}, 
such that, for each $i\in\{1,\cdots,n\}$, the number of birth/death moves  is equal to the number of saddle moves, 
in the deformation from the $i$th component of $L$ into the $i$th component of $L'$.
In the case of based diagrams, the move of Figure~\ref{BC} is also required, and 
each local move does not include base points.

It is shown that Milnor invariants for diagrams are welded-concordance invariants \cite{C}.

\medskip
\begin{remark}
If $L$ and $L'$ are classical links, then the usual {\em link concordance} implies 
the welded-concordance.
In the definition of the welded-concordance, the condition on the numbers of birth/death and saddle moves
corresponds to the fact that each components of $L$ and $L'$ bounds an annulus in the 4-space.
In contrast to the usual link concordance, all welded knots are welded-concordant to the trivial knot \cite{Gaudreau}. 
\end{remark}

\medskip
\subsection{Ascending presentation and Milnor invariants}

In this subsection, we give a sketch of proof of Theorem \ref{Boris} using 
a special $W$-tree presentation called {\em ascending} presentation.

Let $(V, T)$ be a $W$-tree presentation of a based diagram $(L,\mathbf{p})$. 
Note that $V$ is also a based diagram. 
The $W$-tree presentation $(V,T)$ is  {\em ascending}\footnote{This notion was first defined in \cite{ABMW}, in the context of Gauss diagrams. }  if, when running around each component of $V$ from the base point along the orientation, all tails of $T$ appear before all heads of $T$.
A based diagram is  {\em ascendable} if it has an ascending presentation.

\medskip
\subsubsection{Peripheral systems of ascending presentations}

Let $(L,\mathbf{p})$ be an ascendable diagram with an ascending $W$-tree presentation $(V,T)$.
Since $V\cup T$ contains no classical crossings, by (after doing expansions of $W$-trees) arrow moves 
we may assume that $V$ is a based diagram $(O,\mathbf{p})$, where 
$O$ is a trivial diagram, i.e., it contains no crossings. 
In the following, we suppose that $V=(O,\mathbf{p})$.

By the expansions of $W$-trees,  we obtain a $W$-arrow presentation $((O,\mathbf{p}),A)$ from 
the ascending presentation $((O,\mathbf{p}),T)$, where 
for each $W$-arrow in $A$, a small segment adjacent to the head is always assumed to be on the right side with respect to the orientation of the component of $O$ that contains the head.
Furthermore, we may assume that each $W$-arrow contains at most one twist. 
We remark that $((O,\mathbf{p}),A)$ is ascending, since $((O,\mathbf{p}),T)$ is ascending. 

For $((O,\mathbf{p}),A)$, we assign a label by arc $\alpha_i$ to the $W$-arrow whose tail is in the $i$-th component of $O$.
Let $F$ be the free group generated by $\{\alpha_1,...,\alpha_n\}$, 
let $u_{i1},...,u_{ir_i}$ be the labels of $W$-arrows that appear in order,  
when traveling around the $i$th component of $V$ from the base point along the orientation.

Then we define 
\[l^A_i:=u_{i1}^{\e(i1)}\cdot\cdots\cdot u_{ir_i}^{\e(ir_i)}\in F, \]
where $\varepsilon(ij)=-1$ if the $W$-arrow corresponding to $u_{ij}$ contains a twist, and 
 $\varepsilon(ij)=+1$ otherwise.
By the definition of surgery along $W$-arrows, 
it can be seen that  $l^A_i$ is  the sequence of signed labels of the over crossings that appear when running around the $i$th component of $O_A$ from the base point.
Moreover, if necessary, by using AR8 move with preserving the ascending presentation, $l^A_i$ can be regarded 
as an $i$th longitude $\lambda_i^L$.
Since $l_i^A$ is an element of the free group $F$, 
we have the following proposition. 

\medskip
\begin{prop}(\cite{AM})\label{ascending}
Let $((O,\mathbf{p}),A)$ and $((O,\mathbf{p}),A')$ be 
ascending presentations for  ascendable based diagrams 
$(L,\mathbf{p})$ and $(L',\mathbf{p}')$ respectively. 
If $l_i^A=l_i^{A'} (\in F)$ for each $i$, then 
$((O,\mathbf{p}),A)$ and $((O,\mathbf{p}),A')$ are deformed into each other by 
a combination of AR7 and AR9, and hence, 
$(L,\mathbf{p})$ and $(L',\mathbf{p}')$ are equivalent. 
\end{prop}

\medskip
On the other hand, by Proposition~\ref{Magnus}~(2), we have the following. 

\medskip
\begin{prop}(\cite{Col})\label{ascending2}
For any sequence $I$ with length at most $q$, $\mu_{(L,\mathbf{p})}(I)=\mu_{(L',\mathbf{p}')}(I)$ 
if and only if for each $i$, 
$\lambda^L_i\equiv \lambda^{L'}_i \pmod{\Gamma_q F}$.\footnote{This proposition does not 
require the assumption of \lq ascending\rq.}
\end{prop}

\medskip
\subsubsection{Sketch of the proof of Theorem~\ref{Boris}}

Since the \lq only if\rq\  part of  Theorem \ref{Boris} is not difficult to show,
we admit that  the  \lq only if\rq\ part holds and show the \lq if\rq\ part.
Also, since the case $k=1$ is obvious, we suppose that $k\geq 2$.

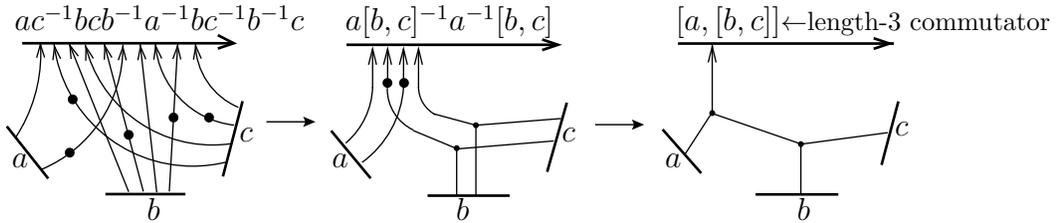
\begin{figure}[!h]
  \begin{center}
    \begin{overpic}[width=12cm]{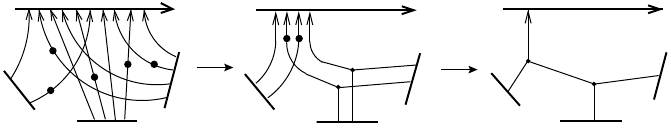}
      \put(4,13){$a$}
      \put(55,-7){$b$}
      \put(90,23){$c$}
      \put(5,65){$a
      {c^{-1}bcb^{-1}}a^{-1}
      {bc^{-1}b^{-1}c}$}
      \put(125,13){$a$}
      \put(172,-7){$b$}
      \put(212,23){$c$}
      \put(130,65){$a
      {[b,c]^{-1}}a^{-1}
      {[b,c]}$}
      \put(251,13){$a$}
      \put(301,-7){$b$}
      \put(338,25){$c$}
      \put(255,65){$[a,[b,c]]$\footnotesize{$\leftarrow$length-$3$ commutator}}
    \end{overpic}
  \end{center}
\caption{Inverse of expansion (where $[x,y]:=xy^{-1}x^{-1}y$)}  \label{reextension}
\end{figure}

\medskip
\begin{enumerate}
\item[{(Step 1)}] 
We deform by $W_k$-concordance the $W$-tree presentations of $(L,\mathbf{p})$ and $(L',\mathbf{p}')$ 
into ascending presentations $((O,\mathbf{p}),T)$ and $((O,\mathbf{p}),T')$, respectively.  
This can be done by using the move (HTE) and welded-concordance. 
Doing the expansions of $T$ and $T'$, we have ascending $W$-arrow presentations 
$((O,\mathbf{p}),A)$ and $((O,\mathbf{p}),A')$. 

\item[{(Step 2)}] 
Since $\mu_{(L,\mathbf{p})}(I)=\mu_{(L',\mathbf{p'})}(I)$ for any sequence 
$I$ with length at most $k$, by the \lq only if\rq\ part of Theorem~\ref{Boris}, 
$\mu_{(O_A,\mathbf{p})}(I)=\mu_{(O_{A'},\mathbf{p'})}(I)$. 
Therefore, by Proposition~\ref{ascending2}, for each $i$, we have  
\[l^{A}_i\equiv l^{{A'}}_i\pmod{\Gamma_{k}F}.\]
In other words, there is an element $g_i$ in $\Gamma_{k}F$ such that $l^{{A'}}_i=l^{A}_i g_i\in F$. 
By Proposition~\ref{ascending}, there are $W$-arrows $B$
such that $O_{A'}=O_{A\cup B}$ and $l_i^{A'}=l_i^Al^B_i$.

\item[{(Step 3)}] 
From an algebraic point of view, we see that $g_i$ is a product of commutators $g_{i1},...,g_{it_i}$ 
consisting of $\{\alpha_1^{\pm1},...,\alpha_n^{\pm1}\}$ with length at least $k$.
And 
by applying the inverse of expansion (Figure~\ref{reextension})  
to the $W$-arrows corresponding to each $g_{is}$, we obtain a $W$-tree of degree at least $k$. 
Hence $O_{A'}=O_{A\cup B}$ and $O_A$ are $W_k$-equivalent.  
\end{enumerate}

\medskip
\subsection{Self $W_k$-concordance and Minor invariants}

The content of this section is due to  \cite{AMY1}.
The {\em self $W_k$-concordance} is an equivalence relation on diagrams obtained by combining 
self $W_k$-equivalence and welded-concordance.
Here, we introduce the result that the self $W_k$-concordance classification is given by Milnor invariants.

\medskip
\subsubsection{$k$-reduced quotient groups and Milnor invarints}


In this section, we define the $k$-reduced quotient group of a based diagram $(L,\mathbf{p})$, 
and explain the relation between this group and Milnor invariants $\mu_{(L,\mathbf{p})}(I)~(r(I)\leq k)$, where 
\[r(I):=\max\{r_i~|~\text{$i$ appears $r_i$ times in $I$}~(i=1,...,n)\}.\] 
The arguments in this section are simply rewritings of those in  \cite{AMY1} for welded 
string links into those for based diagrams.

For the peripheral system $(G(L,\mathbf{p}),\{a_{i0},\lambda_i\}_i)$ of $(L,\mathbf{p})$, 
let $G_i(L,\mathbf{p})$ be the normal closure of $ G(L,\mathbf{p})$ that contains $a_{i0}$.
We call the quotient group 
\[R_k(L,\mathbf{p}):= G(L,\mathbf{p}) / \Gamma_{k+1} G_1(L,\mathbf{p})\cdots 
\Gamma_{k+1} G_n(L,\mathbf{p})\] 
the {\em $k$-reduced group of $(L,\mathbf{p})$} or the {\em $k$-reduced quotient of $G(L,\mathbf{p})$}.

When $k=1$, it is the {\em link group} defined by Milnor \cite{M1}.
Since $G(L,\mathbf{p})$ is also sometimes called link group, we use a different name here.
Habegger and Lin \cite{HL} call it the {\em reduced group}.

It is known that $R_k(L,\mathbf{p})$ is a nilpotent group for any positive integer $k$, 
and 
\[\Gamma_{kn+1} G(L,\mathbf{p})
\subset \Gamma_{k+1} G_1(L,\mathbf{p})\cdots  \Gamma_{k+1} G_n(L,\mathbf{p}).\]
Therefore we have the following surjection.
\[
N_{kn+1}(L,\mathbf{p})
=G(L,\mathbf{p})/\Gamma_{kn+1}G(L,\mathbf{p})
\twoheadrightarrow 
R_k(L,\mathbf{p}).
\]
By composing respectively this with the two maps 
\[\rho_{kn+1}:G(L,\mathbf{p})\rightarrow N_{kn+1}(L,\mathbf{p}),
~~~~\phi_{kn+1}:A(L,\mathbf{p})\rightarrow N_{kn+1}(L,\mathbf{p})\]
defined in Sections~\ref{peri-diagram} and \ref{color-diagram}, we have the two maps 
\[ \rho_{R_k}:G(L,\mathbf{p})\rightarrow R_k(L,\mathbf{p}),~~~~
\phi_{R_k}:A(L,\mathbf{p})\rightarrow R_k(L,\mathbf{p}).\]

Since $N_{kn+1}(L,\mathbf{p})$ is generated by $\rho_{kn+1}(a_{i0})~(i=1,...,n)$, 
$R_k(L,\mathbf{p})$ is generated by $\rho_{R_k}(a_{i0})=\alpha_i~(i=1,...,n)$. 
Hence $\rho_{R_k} (\lambda_k)$ is written as a word of $\alpha_1,...,\alpha_n$. 

Let $F$ be the free group generated by $\alpha_1,...,\alpha_n$, 
and let $N_i$ be the normal closure containing $\alpha_i$ for each $i=1,...,n$. 
Similar to Theorem~\ref{CM}~(1), the following theorem holds.

\medskip
\begin{thm}\label{AMY0}~
The Chen-Milnor map $\eta_{kn+1}$ induces the following isomorphism 
\[R_k(L,\mathbf{p})\cong\langle\alpha_1,\ldots,\alpha_n~|~ 
\Gamma_{k+1}N_1\cdots\Gamma_{k+1}N_n\rangle~(=F/\Gamma_{k+1} N_1\cdots\Gamma_{k+1} N_n).\]
\end{thm}

\medskip
As in Section~\ref{color-diagram}, by the theorem above, 
$\phi_{R_k}$ is a $(F/\Gamma_{k+1} N_1\cdots\Gamma_{k+1} N_n)$-coloring.

\medskip
\subsubsection{Self $W_k$-concordance classification}

As in Section~\ref{aut-nil-diagram}, we can define the automorphism 
\[\varphi_{R_k}:R_kF\longrightarrow R_kF,~\alpha_k\longmapsto 
\rho_{R_k}(\lambda_k)^{-1}\alpha_k\rho_{R_k}(\lambda_k)
\]
of $R_kF=F/\Gamma_{k+1} N_1\cdots\Gamma_{k+1} N_n$.

The following theorem gives the self $W_k$-concordance classification of based diagrams.

\medskip
\begin{thm}\label{AMY}
Let $(L,\mathbf{p})$ and $(L',\mathbf{p}')$ be based diagrams.
Then for any positive integer $k$, the following (1),(2) and (3) are mutually equivalent.
\begin{enumerate}
\item For any sequence $I$ with $r(I)\leq k$, $\mu_{(L,\mathbf{p})}(I)=
\mu_{(L',\mathbf{p}')}(I)$.
\item $\varphi_{R_k}(L,\mathbf{p})=\varphi_{R_k}(L',\mathbf{p}')\in 
\mathrm{Aut}_{\mathrm{c}}(R_kF)$.
\item $(L,\mathbf{p})$ and $(L',\mathbf{p}')$ are self $W_k$-concordant.
\end{enumerate}
Here 
$\mathrm{Aut}_{\mathrm{c}}(R_kF)$ is the set of automorphisms 
of $R_kF$ that act by conjugation on each generator. 
\end{thm}

\medskip
\begin{remark}
(1)~In \cite{AMY1}, the statement of Theorem~\ref{AMY} is given for welded string links 
rather than based diagrams. These statements are essentially the same.
The self $W_1$-concordant classes coincide with the self $W_1$-equivalence classes.
The self $W_1$-equivalence classification of welded string links is given in  \cite{ABMW}.\\
(2)~The set of self $W_k$-concordance classes $\mathrm{wSL}/{(\text{self }W_k+\mathrm{c})}$ 
of welded $n$-string links forms a group under the product of welded string links. 
Moreover by using Theorem~\ref{AMY} (for the welded string link version), 
we can show the following isomorphism
\[\mathrm{wSL}/{(\text{self }W_k+\mathrm{c})} \cong \mathrm{Aut}_{\mathrm{c}}(R_k F).\]
\end{remark}

\medskip
\begin{remark}
It is shown by T. Fleming and the author \cite{FY} that 
for any sequence $I$ with $r(I)\leq k$, Milnor invariant $\mu(I)$  
is a self $C_k$-equivalence invariant.
The self $C_1$-equivalence is  link-homotopy, and 
the self $C_2$-equivalence is often called the {\em self $\Delta$-equivalence}. 
For the case $k=1,2$, it is shown in \cite{M1} and \cite{yasuhara} respectively 
that a classical link $L$ is self $C_k$-equivalent to a trivial 
if and only if $\omu_L(I)=0$ for any sequence $I$ with $r(I)\leq k$.
The self $C_k$-concordance classification for classical string links 
is given in \cite{HL} for $k=1$ and in \cite{yasuharaAGT} for $k=2$. 
\end{remark}

\medskip
\subsubsection{Sketch of the proof of Theorem~\ref{AMY}}

The key for the proof of Theorem~\ref{AMY} is ascending presentations, as in the proof of Theorem~\ref{Boris}. 
We need two additional theorems below.

A $W$-tree $T$ for a diagram $L$ is a {\em $W^{(k)}$-tree} if 
at least $k$ ends of $T$ belong to the same component of $L$.
The {\em $W^{(k)}$-equivalence} is an equivalence relation on diagrams generated by 
surgery along $W^{(k)}$-trees, and 
the {\em $W^{(k)}$-concordance} is a combination of the $W^{(k)}$-equivalence and welded-concordance.

\medskip
\begin{thm}\label{AMY1}
Let $k$ be a positive integer. 
Two diagrams are $W^{(k+1)}$-equivalent if and only if they are self $W_k$-equivalent.
\end{thm}

\medskip
\begin{thm}\label{AMY2}
Let $\lambda_i^L$ and $\lambda_i^{L'}$ be $i$th longitudes 
of based diagrams $(L,\mathbf{p})$ and $(L',\mathbf{p'})$, respectively.  
Then $\mu_{(L,\mathbf{p})}(I)=\mu_{(L',\mathbf{p})}(I)$ for 
any sequence $I$ with $r(I)\leq k$ if and only if 
$\rho_{R_k}(\lambda^L_i)\equiv\rho_{R_k}(\lambda^{L'}_i)
 \pmod{ \Gamma_{k+1}N_1\cdots \Gamma_{k}N_i\cdots \Gamma_{k+1}N_n}$.
\end{thm}

\medskip
Since the proof of $(1)\Rightarrow(3)$ is the most difficult,   
we admit $(3)\Rightarrow(1)$ and $(1)\Leftrightarrow (2)$ (and Theorems~\ref{AMY1} and \ref{AMY2}) 
and only give a sketch of proof of $(1)\Rightarrow(3)$.

\medskip
\begin{enumerate}
\item[{(Step 1)}] 
We deform by $W_{kn}$-concordance the $W$-tree presentations of $(L,\mathbf{p})$ and $(L',\mathbf{p}')$ 
into ascending presentations $((O,\mathbf{p}),T)$ and $((O,\mathbf{p}),T')$, respectively.  
Doing expansions of $T$ and $T'$, we have ascending $W$-arrow presentations 
$((O,\mathbf{p}),A)$ and $((O,\mathbf{p}),A')$. 
Since a $W_{kn}$-tree has $kn+1$ ends, 
it is a $W^{(k+1)}$-tree.  Hence, by Theorem~\ref{AMY1}, 
$(L,\mathbf{p})$ and $(L',\mathbf{p}')$ are self $W_k$-equivalent to 
$((O,\mathbf{p}),A)$ and $((O,\mathbf{p}),A')$ respectively. 

\item[{(Step 2)}] Since $\mu_{(L,\mathbf{p})}(I)=\mu_{(L',\mathbf{p'})}(I)$ for any sequence 
$I$ with $r(I)\leq k$, by $(3)\Rightarrow(1)$ of Theorem~\ref{AMY}, 
$\mu_{(O_A,\mathbf{p})}(I)=\mu_{(O_{A'},\mathbf{p'})}(I)$. 
Therefore, by Theorem~\ref{AMY2}, for each $i$, we have  
\[l^{A}_i\equiv l^{{A'}}_i\pmod{ \Gamma_{k+1}N_1\cdots \Gamma_{k}N_i\cdots \Gamma_{k+1}N_n}.\]
Hence, there are elements $g_{ij}\in \Gamma_{k+1}N_j ~(j\neq i),~g_{ii}\in\Gamma_{k}N_i$ 
 such that $l^{{A'}}_i=l^{A}_i g_{i1}\cdots g_{in}\in F$. 
By Proposition~\ref{ascending}, there are $W$-arrows $B$
such that $O_{A'}=O_{A\cup B}$ and $l_i^{A'}=l_i^Al_i^B$. 

\item[{(Step 3)}] 
From an algebraic point of view, we have that
\begin{enumerate}
\item if $j\neq i$, then 
 $g_{ij}\in \Gamma_{k+1}N_j$ is a product of commutators $g_{ij}^{1},...,g_{ij}^{s_j}$ 
 with length at least $k$ where $\alpha_j$ appears at least $k+1$ times, and 
\item if $j= i$, $g_{ii}^{k}\in \Gamma_{k}N_i$ is a product of commutators $g_{ii}^{1},...,g_{ii}^{s_i}$
 with length at least $k$ where $\alpha_i$ appears at least $k$ times.
\end{enumerate}

\medskip\noindent
In both cases, 
by applying the inverse of expansion (Figure~\ref{reextension})  
to the $W$-arrows corresponding to each $g_{is}$, we obtain a $W$-tree 
such that at least $k+1$ ends belong to the same component of $O$. 
Hence, by Theorem~\ref{AMY1},
$O_{A'}=O_{A\cup B}$ and $O_A$ are self $W_k$-equivalent.  
\end{enumerate}

\section{Milnor invariants for surface-links}

The content of this and the next chapters is taken from \cite{AMY}. 
Milnor invariants can be defined for {\em $m$-dimensional links}, which are 
$m$-manifolds (smoothly embedded) in $(m+2)$-space \cite{AMY}. 
For simplicity, we explain the case $m=2$.
A 2-dimensional link often means spheres in $4$-space. 
Milnor invariants introduced here vanish if the fundamental groups of the objects themelves are trivial. 
Here we consider oriented surfaces, not only spheres, in $4$-space, and call them {\em surface-links}.
It is possible to consider surfaces with boudary, mimicking string links, 
but to avoid complications we consider only closed surfaces.

\medskip
\subsection{Surface-link diagrams}

\begin{figure}[!h]
\begin{center}
\includegraphics[width=.33\linewidth]{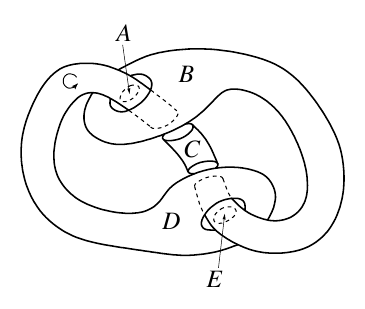}
\caption{A surface-link (surface-knot) diagram}\label{Sphere}
\end{center}
\end{figure}

As a classical link in $3$-space can be drawn as a diagram in the plane,
a surface-link in $4$-space can be drawn as a diagram in $3$-space, 
for example, see Figure \ref{Sphere}.  
(The arrowed circle $\circlearrowleft$ in the figure represents the {\em surface orientation}.
When turning a right hand screw in the direction of $\circlearrowleft$, 
by specifying   \lq front\rq\ and \lq back\rq\  of the surface so that the screw moves back to front, 
the \lq orientation\rq\ of the surface can also be seen as a \lq specifying front/back\rq.)
We note that, in the case of surface-links, the intersections of surfaces appear as curves, 
and over/under informations are added according to the arrangements in $4$-space 
by cutting along \lq under-intersections\rq.  

\begin{figure}[!h]
\begin{center}
\includegraphics[width=.25\linewidth]{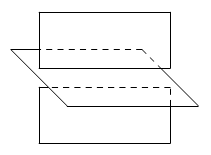}
\hspace{1cm}
\includegraphics[width=.17\linewidth]{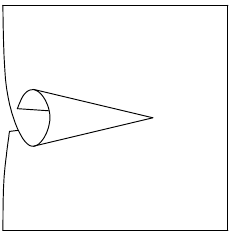}
\hspace{1cm}
\includegraphics[width=.25\linewidth]{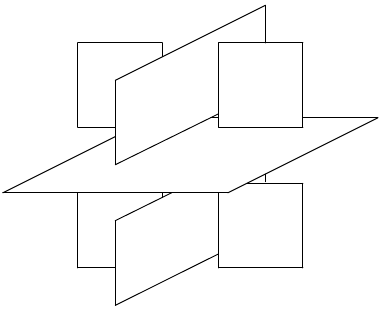}
   \vspace{.5em}
\caption{Intersections of diagrams}\label{sing}
\end{center}
\end{figure}

In Figure~\ref{Sphere}, only \lq simple\rq\ intersections of the surface appear as in the left side of Figure~\ref{sing}.
But in general, 
as in the center of Figure~\ref{sing}, an intersection containing a point called a {\em branch point} may appear, 
and as shown on the right side of Figure~\ref{sing}, three surfaces may intersect, forming a {\em triple point}.
On the other hand, it is known that by using these three types of intersections, 
any surface-link can be drawn in 3-space.\footnote{There are two types of branch points, positive and 
negative, see Figure~\ref{CutBranch}.}
A diagram drawn in this way is called a {\em diagram of a surface-link}. 

As there is Reidemeister's theorem for classical links, there is also the following theorem for surface-links.
The local moves on diagrams of surface-links that correspond to  Reidemeister moves of classical link diagrams 
are called  {\em Roseman moves}.

\medskip
\begin{thm}(Roseman's Theorem \cite{Rose,CS,CS2})
Two surface-links are equivalent (i.e., there is a continuous deformation of 4-space between them)
if and only if their diagrams are deformed into each other 
by a combination of continuous deformations of 3-space and Roseman moves. 
\end{thm}

\medskip
In general, by allowing the three kinds of intersections in Figure~\ref{sing},  
($n$-component) surfaces drawn in 3-space are called ($n$-component) {\em surface-link diagrams}.
When $n=1$, we often call them {\em surface-knot diagrams}. 
Figure~\ref{Sphere} shows a surface-knot diagram.
A surface-link diagram is {\em trivial} if it has no intersections.

A surface-link diagram $L$ can be separated into several {\em regions} by under-intersections. 
For example, the diagram in Figure~\ref{Sphere} is separated into 5 regions $A,B,C,D,E$. 
For an $n$-component surface-link diagram $L$ and for each $i~(i=1,...,n)$, 
we choose a point $p_i$ on a region of the $i$th component. 
Then we call the pair $(L,\mathbf{p})$ of $L$ and $\mathbf{p}=(p_1,...,p_n)$ a
{\em based surface-link diagram}.
Two based surface-link diagrams are {\em equivalent} 
if one is deformed into the other 
by a combination of continuous deformations of 3-space and Roseman moves 
that do not contain base points. 

\medskip
\subsection{Groups of based surface-link diagrams}

\begin{figure}[!h]
  \begin{center}
    \begin{overpic}[width=4.2cm]{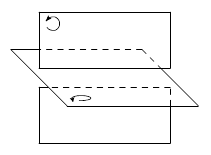}
      \put(48,65){$x$}
      \put(48,15){$y$}
      \put(98.5,29){$z$}
    \end{overpic}
  \end{center}
  \caption{Relation~${z}xz^{-1}{y}^{-1}~(x,y,z
  \in\{a_{ij}\}_{ij})$}
  \label{relation}
\end{figure}

Let $(L,\mathbf{p})$ be a based surface-link diagram. 
We define the group $G(L,\mathbf{p})$ of $(L,\mathbf{p})$ as follows.
For the $i$th component of $L$, let $a_{i0}$ be the region that contains the base point $p_i$, 
and let $a_{i1},...,a_{i{r_i}}$ be the other regions, where  
$r_i+1$ is the number of regions in $K_i$. 
(Unlike classical link diagrams, there are no rules for labeling the regions $a_{ij}~(j>0)$, 
and only the region $a_{i0}$ contains the base point $p_i$.) 
Let $\widetilde{F}$ be the free group generated by $\{a_{ij}\}_{i,j}$.  
For each simple intersection of $L$, we consider a relation
$zxz^{-1}{y}^{-1}$ as in Figure~\ref{relation}. 
The {\em group $G(L,\mathbf{p})$ of $(L,\mathbf{p})$} is the quotient group of the 
free group $\widetilde{F}$ with a generating set $\{a_{ij}\}_{i,j}$ modulo 
these relations.  
It is known that $G(L,\mathbf{p})$ is isomorphic to the fundamental group of 
the complement of the surface-link whose diagram is $L$. 

\medskip
\subsection{Nets of diagrams}

In this section,  we intuitively  explain {\em cut-diagrams} of  
classical links and of surface-links. 
As we see below, cut-diagrams are \lq nets\rq\ of diagrams. 

\medskip
\subsubsection{Cut-diagrams of classical link diagrams}

\begin{figure}[!h]
\begin{center}
\includegraphics[width=.20\linewidth]{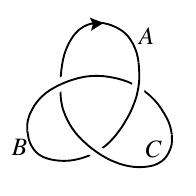}\hspace{1cm}
\includegraphics[width=.22\linewidth]{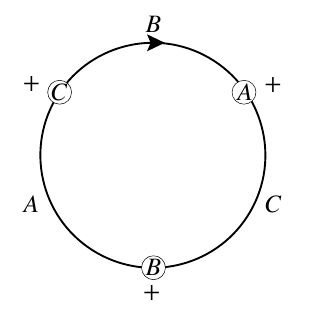}
\caption{A link (knot) diagram and its cut-diagram}\label{Trefoil}
\end{center}
\end{figure}

First, we explain cut diagrams of classical link diagrams.
The diagram on the left side of Figure~\ref{Trefoil} is obtained from the circle on the right side 
by \lq gluing\rq\ points in the arcs $A, B, C$ to the points $\text{\textcircled{\raisebox{-1.1pt}{$A$}}}$, 
$\text{\textcircled{\raisebox{-1.1pt}{$B$}}}$, $\text{\textcircled{\raisebox{-1.1pt}{$C$}}}$ respectively 
so that $\text{\textcircled{\raisebox{-1.1pt}{$A$}}}$, 
$\text{\textcircled{\raisebox{-1.1pt}{$B$}}}$, $\text{\textcircled{\raisebox{-1.1pt}{$C$}}}$ become 
under crossings, and the signs of crossings coincide with the signs of the points. 
From this observation, it can be seen that the circle  on the right side of Figure~\ref{Trefoil}
is a net of the diagram on the left side.
We call such nets {\em cut-diagrams of link diagrams}. 
That is, a cut diagram of a link diagram is a set of circles with signed points, that correspond 
to under crossings and are   labeled with arcs containing their over crossings.

\medskip
\begin{remark}\label{uniqueness}
For the cut-diagram in Figure~\ref{Trefoil}, 
let $T$ be a set of three $W$-arrows $\alpha_X~(X=A,B,C)$ such that 
the tail of $\alpha_X$ belongs to the arc $X$ and the head is on the point $\text{\textcircled{\raisebox{-1.1pt}{$X$}}}$.  
And let $(O,T)$ be an arrow presentation consisting of the trivial knot diagram 
$O$ and $T$, where 
for each $W$-arrow in $\alpha_X$, a small segment adjacent to the head is always assumed to be 
on the right side with respect to the orientation of $O$, and 
the number of twists is 0 if the sign of the point $X$ is positive, or
1 otherwise.
Then $(O,T)$ is an arrow presentation of the knot diagram on the left side of Figure~\ref{Trefoil}. 
(In general, since $T$ may contain more than one $W$-arrows whose tails belong the same arc, 
$(O,T)$ is not uniquely determined.)
Likewise, 
from a cut-diagram of any link diagram, we obtain an arrow presentation on a trivial link diagram $O$.  
For a link diagram $L$, 
let $(O,T)$ be an arrow presentation obtained from a cut-diagram of $L$,
 and let  $(O,T')$ be an arrow presentation of $L$ obtained by the operations in Figure~\ref{Aprst}.  
Although the heads of $T$ and $T'$ coincide, 
their tails may not coincide. 
On the other hand, 
 thanks to AR7 (VR1-VR3 and AR6 if necessary),  
$(O,T)$ is equivalent to $(O,T')$, and hence it is an arrow presentation of $L$. 
It follows that two link diagrams that have a common cut-diagram
are equivalent as welded links. 
\end{remark}

\medskip
\subsubsection{Cut-diagrams of surface-link diagrams}

\begin{figure}[!h]
\begin{center}
\includegraphics[width=.3\linewidth]{Sphere.pdf}
\includegraphics[width=.2\linewidth]{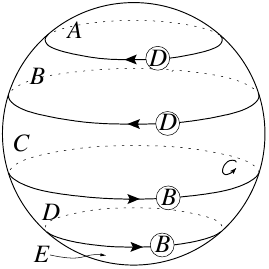}
  \vspace{.5em}
\caption{A surface-link diagram (spherical-link diagram) and its cut-diagram}
\label{CutSphere}
\end{center}
\end{figure}

By applying the idea of \lq net\rq\  to surface-link diagrams,  we obtain 
{\em cut-diagrams of surface-link diagrams}.
The figure on the right side of Figure~\ref{CutSphere} is a cut-diagram of the surface-link diagram on the left side of 
Figure~\ref{CutSphere}. 
In contrast to classical link diagrams, for surface-link diagrams, 
since the gluing parts are curves, the gluing is specified using oriented, labeled curves on surfaces  
that correspond to under-intersections.\footnote{These gluing parts without labels 
are also called lower decker sets \cite{CS2}.}
 
For example, in the cut-diagram on the right side of Figure~\ref{CutSphere},  
the curve with label $D$ between the regions $A$ and $B$ intersects the region $D$ and is oriented so that, 
in the gluing result (the left side of Figure~\ref{CutSphere}),  
\[(\text{direction of back to front of $D$},
\text{direction of back to front of $A$(and $B$)},
\text{orientation of $\text{\textcircled{\raisebox{-1.1pt}{$D$}}}$})\]
is equal to the right-handed spatial coordinate system 
\[(\text{direction of $x$-axis},\text{direction of $y$-axis},\text{direction of $z$-axis}).\]

\begin{figure}[!h]
\begin{center}
\begin{tabular}{rcl}
\begin{minipage}[t]{.15\hsize}
\centering\includegraphics[width=1\linewidth]{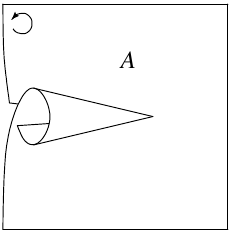}
\end{minipage}
&
$\begin{array}{c}
\stackrel{\text{net}}{\longrightarrow}\\[1.7cm]
~
\end{array}$
&
\begin{minipage}[t]{.15\hsize}
\centering\includegraphics[width=1\linewidth]{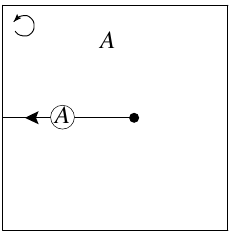}
\end{minipage}
\end{tabular}\hspace{1cm}
\begin{tabular}{rcl}
\begin{minipage}[t]{.15\hsize}
\centering\includegraphics[width=1\linewidth]{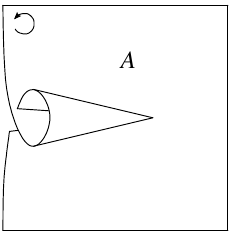}
\end{minipage}
&
$\begin{array}{c}
\stackrel{\text{net}}{\longrightarrow}\\[1.7cm]
~
\end{array}$
&
\begin{minipage}[t]{.15\hsize}
\centering\includegraphics[width=1\linewidth]{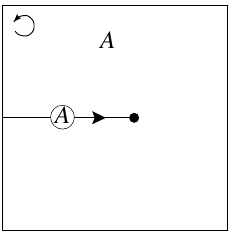}
\end{minipage}
\end{tabular}
\vspace*{-1cm}
\caption{Cut-diagrams (local nets) around branch points}
\label{CutBranch}
\end{center}
\end{figure}

\begin{figure}[!h]
\begin{center}
\begin{tabular}{rcl}
\begin{minipage}[t]{.25\hsize}
\centering
\includegraphics[width=1\linewidth]{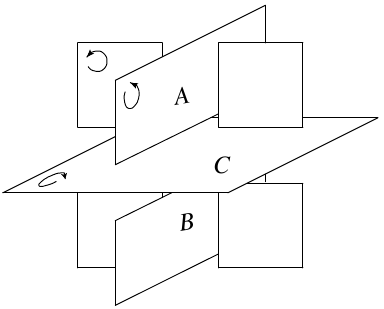}
\end{minipage}
&
$\begin{array}{c}
\stackrel{\text{net}}{\longrightarrow}\\[2.5cm]
~
\end{array}$
&
\begin{minipage}[t]{.6\hsize}
\centering
\includegraphics[width=1\linewidth]{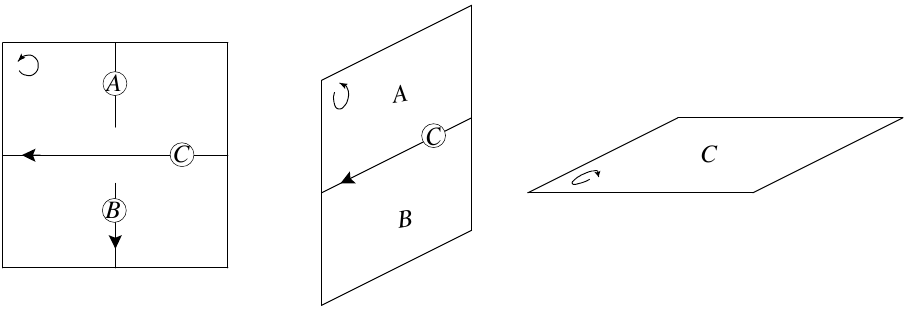}
\end{minipage}
\end{tabular}\\[-1.5cm]
\begin{tabular}{rcl}
\begin{minipage}[t]{.25\hsize}
\centering
\includegraphics[width=1\linewidth]{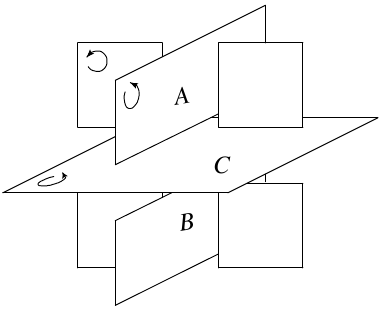}
\end{minipage}
&
$\begin{array}{c}
\stackrel{\text{net}}{\longrightarrow}\\[2.5cm]
~
\end{array}$
&
\begin{minipage}[t]{.6\hsize}
\centering
\includegraphics[width=1\linewidth]{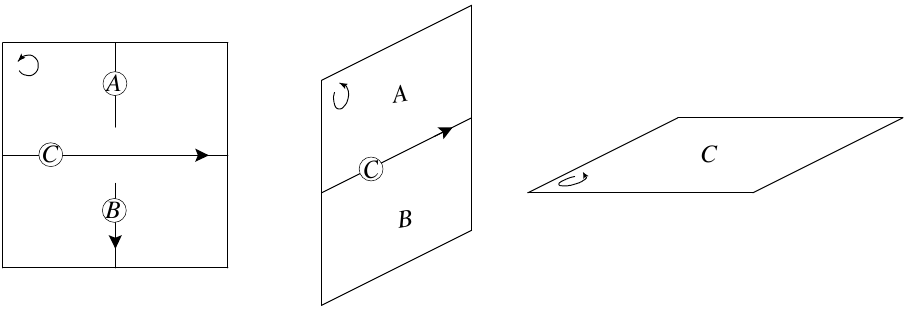}
\end{minipage}
\end{tabular}
\vspace*{-1.5cm}
\caption{Cut-diagrams (local nets) around triple points}\label{CutTriple}
\end{center}
\end{figure}

In general, surface-link diagrams may contain branch points and triple points. 
Then cut-diagrams have \lq local nets\rq\  as illustrated in 
Figures~\ref{CutBranch} (resp. \ref{CutTriple}), which contain {\em end points} corresponding to 
the branch points (resp. classical crossings corresponding to the triple points).

\medskip
\subsection{Cut-diagrams}

In the previous section, we defined cut-diagrams from (surface-)link diagrams. 
But, in this section, we define cut diagrams independently of (surface-)link diagrams.

\medskip
\subsubsection{$1$-dimensional cut-diagrams}
A cut diagram of an $n$-component link diagram can be seen as a set  $S$ 
of $n$ circles ($1$-dim) with signed points $P$ ($0$-dim) arranged on $S$, where each point in $P$ is labeled by an arc of $S-P$.
Here, instead of having cut-diagrams from link diagrams,
we consider cut-diagrams independently of link diagrams, that are called  {\em $1$-dimensional cut-diagrams}. 
That is, an $n$-component 1-dimensional cut-diagram is  
a triple $(S,P,f)$  consisting of a set $S$ of $n$ circles, a set of signed points $P$ on $S$, and 
a map $f$ from a set of arcs of $S-P$ to $P$. 
For a cut-diagram $C(=(S,P,f))$, we can build a diagram from $C$ by gluing it. 
First, we locally glue only the crossings according to the correspondence of $f$ to create classical crossings, and then connect these crossings with corresponding arcs to complete the gluing. 
If $C$ is cut-diagram of a certain link diagram, then we can connect crossings so that 
the result is equal to the original link diagram.  
But, in general, when connecting crossings on a plane, new crossings may be necessary. 
In such cases, a virtual link diagram can be obtained by making the new crossings virtual crossings.
Therefore, for any cut-diagram $C$, by gluing it, we have a virtual link diagram $L_C$. 
Although $L_C$ might not be uniquely determined, by the observation in Remark~\ref{uniqueness}, 
it is unique up to welded moves. 
Since we can obtain a cut-diagram from any virtual link-diagram, we have the following 
sequence of surjections. 
\[\{\text{virtual link diagrams}\}\twoheadrightarrow\{\text{cut-diagrams}\}
\twoheadrightarrow\{\text{welded links}\}.\]
While the surjection from cut-diagrams to welded links is not injective, 
it induces a bijection from a quotient of cut-diagrams modulo certain moves. 
Here these moves are just the direct translations of welded moves. 
In this sense, 1-dimensional cut diagrams can be seen as welded link diagrams. 
\medskip
\subsubsection{$2$-dimensional cut-diagrams}

A cut-diagram of an $n$-component surface-link diagram can be seen as 
a set $\Sigma$ of $n$ closed surfaces ($2$-dim) with oriented 1-dimensional diagram\footnote{This is almost 
a link diagram but it may contain end points that correspond to branch points}$D$ (1-dim), 
where each arc of $D$ is labeled by a region of $\Sigma-D$.
Therefore we define  {\em $2$-dimensional cut-diagrams} independently of surface-link diagrams 
as a triple $(\Sigma,D,f)$ consisting of 
a set $\Sigma$ of $n$ closed surfaces, a $1$-dimensional diagram $D$ on $\Sigma$, and  
a map $f$ from the set of arcs of $D$ to the set of regions of $\Sigma-D$. 
Here, the map $f$ is defined so that \lq local gluing can be performed\rq\  
around each branch point and  each triple point.\footnote{For a detailed definition, see \cite[Subsection 1.2.1]{AMY}.} 
For an arc $\alpha$ of $D$, $f(\alpha)$ is called the {\em label of $\alpha$}. 
By putting a circle $\text{\textcircled{\raisebox{-1.1pt}{$A$}}}$ on 
$\alpha$, we mean that the label $f(\alpha)$ of $\alpha$ is $A$.

Since $1$-dimensional cut-diagrams can be regarded as welded links, $2$-dimensional 
cut-diagrams can be considered as 2-dimensional generalization of welded links.

\medskip
\subsection{$2$-dimensional cut-diagrams and groups}

One important property of cut-diagrams is that they contain gluing information around the intersections. 
From this information, the group of the original surface-link diagram can be computed from the cut-diagram.

In the following, we define peripheral systems of (based) 2-dimensional cut-diagrams and 
their Milnor invariants, which induce  Milnor invariants for surface-links.
From now on, in this chapter, we only treat 2-dimensional cut-diagrams, 
unless otherwise specified, cut-diagrams always mean $2$-dimensional.

\subsubsection{Peripheral systems of based cut-diagrams}

For an $n$-component cut-diagram $C=(\Sigma,D,f)$,
let $\Sigma_i$ be the $i$th component of $\Sigma$, and let $D_i=\Sigma_i\cap  D$.
For each $i$, let $p_i$ be a point on a region of $\Sigma_i-D_i$. 
Then the pair $(C,\mathbf{p})$ of $C$ and $\mathbf{p}=(p_1,...,p_n)$ is called 
a {\em based cut-diagram}.

\begin{figure}[!h]
  \begin{center}
        \begin{overpic}[width=2.5cm]{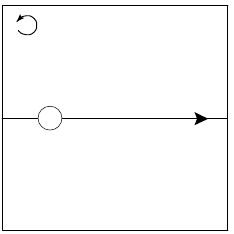}
      \put(33,50){$x$}
      \put(33,15){$y$}
      \put(13,33.5){$z$}
    \end{overpic}
  \end{center}
  \caption{Relation~$zxz^{-1}y^{-1}~(x,y,z
  \in\{a_{ij}\}_{ij})$}
  \label{relation-cut}
\end{figure}

We define the group $G(C,\mathbf{p})$ of $(C,\mathbf{p})$ as follows.
For $\Sigma_i-D_i$, let $a_{i0}$ be the region that contains the base point $p_i$, 
and let $a_{i1},...,a_{i{r_i}}$ be the other regions, where  
$r_i+1$ is the number of regions in $\Sigma_i-D_i$. 
Let $\widetilde{F}$ be the free group generated by $\{a_{ij}\}_{i,j}$.  
For each arc of $D$, we consider a relation
$zxz^{-1}{y}^{-1}$ as in Figure~\ref{relation-cut}. 
The {\em group $G(C,\mathbf{p})$ of $(C,\mathbf{p})$} is the quotient group of the 
free group $\widetilde{F}$ with generating set $\{a_{ij}\}_{i,j}$ modulo 
these relations.  
If $(C,\mathbf{p})$ is a cut-diagram of a based surface-link diagram $(L,\mathbf{p})$, then 
$G(C,\mathbf{p})$ is equal to the group $G(L,\mathbf{p})$ defined in Section~4.2. 
(Since $G(C,\mathbf{p})$ is independent of the choice of $\mathbf{p}$, 
$G(C):=G(C,\mathbf{p})$ is called the {\em group of the cut-diagram $C$}.)

For an oriented loop $l$ on $\Sigma_i$ with base point $p_i$, we define
an element $w(l)$ in $G(C,\mathbf{p})$ as follows. 
Let $u_1,...,u_{l(i)}$ be the sequence of labels of the arcs in $D_i$ that appear in order,  
when traveling along $l$ from $p_i$ following the orientation. 
While traveling around $l$ on the front side of $\Sigma_i$, we define 
 $\varepsilon(j)=+1$ if $u_{j}$ crosses from left to right,
 $\varepsilon(j)=-1$ if $u_{j}$ crosses from right to left.
Then we set 
\[w(l)=a_{i0}^ {-w_i}u_1^{\varepsilon(1)}u_2^{\varepsilon(2)}\cdots u_{l(i)}^{\e(l(i))}\in \widetilde{F},
\]
where $w_i$ is the sum of signs $\e(s)$ for $u_s\subset D_i$.

\medskip
\begin{lem}\label{homotopy-equi}
If two loops $l$ and $l'$ with base points $p_i$ represent the same element of $\pi_1(\Sigma_i,p_i)$, 
then $w(l)$ and $w(l')$ are equal as elements of 
$G(C,\mathbf{p})$.\footnote{By this lemma, the correspondance $l\mapsto w(l)$ induces a map 
$\pi_1(\Sigma_i,p_i)\longrightarrow G(C,\mathbf{p})$. Moreover we have that 
the map is a homomorphism.}
\end{lem}

\medskip
By this lemma, for a loop $l$ representing an element $\pi_1(\Sigma_i,p_i)$, 
$w(l)$ can be seen as an element of $G(C,\mathbf{p})$. 
For each $i$, we call $a_{i0}$  the {\em $i$th meridian}, and we call 
\[\Lambda_i:=\{w(l)~|~\text{$l$ represents a nontrivial element of $\pi_1(\Sigma_i,p_i)$}\}\]
the {\em $i$th longitude set} of $(C,\mathbf{p})$.
The pair $(G(C,\mathbf{p}),\{{a_{i0}},\Lambda_i\}_i)$  
 is called the {\em peripheral system} of $(C,\mathbf{p})$. 

Let $\{l_{ij}~|~j=1,...,2g_i\}$ be a set of loops on $\Sigma_i$ with base points $p_i$
that represent generators of $\pi_1(\Sigma_i,p_i)$, and set 
$\Lambda^0_i:=\{w(l_{ij})~|~j=1,...,2g_i\}$, 
where $2g_i$ is the minimum number of generators of $\pi_1(\Sigma_i,p_i)$, i.e., 
$g_i$ is the genus of $\Sigma_i$. 
By Lemma~\ref{homotopy-equi}, 
any element in $\Lambda_i$ can be written as a product of elements in $\Lambda_i^0$. 
We call $\Lambda_i^0$ an {\em $i$th longitude system}.

\medskip
\subsection{Milnor invariants of $2$-dimensional cut-diagrams}

In this sction, we define Milnor invariants for cut-diagrams, and 
also Milnor invariants for surface-links.

Let $\widetilde{N}$ be the normal closure of $\widetilde{F}$ 
 that contains the relations of (the group presentation of) $G(C,\mathbf{p})$, and 
 let 
$N_q(C,\mathbf{p})=G(C,\mathbf{p})/\Gamma_qG(C,\mathbf{p})$. 
By composing the following sequence 
\[\widetilde{F}\twoheadrightarrow\widetilde{F}/\widetilde{N}
= G(C,\mathbf{p})\twoheadrightarrow 
N_q(C,\mathbf{p})\]
of maps, we have a map
\[\xi_q:\widetilde{F}\longrightarrow N_q(C,\mathbf{p}).\]
It is known that 
$N_q(C,\mathbf{p})$ is a nilpotent group generated by $\xi_q(a_{i0})=\alpha_i~(i=1,...,n)$. 

In the Magnus expansion of $\xi_q (w(l))$,
let $\mu_{(C,\mathbf{p})}(i_1,...,i_s;l)~(s<q)$ be the coefficient of $X_{i_1}\cdots X_{i_s}$. 
For sequences $Ii$ and $J$ with length at most $q$, we define 
\[m_{(C,\mathbf{p})}(Ii):=\gcd\{\mu_{(C,\mathbf{p})}(I;l_{ij})~|~j=1,...,2g_i\}\]
and
\[\Delta_{(C,\mathbf{p})}(J):=\gcd
\left\{m_{(C,\mathbf{p})}(J')~\vline~
\begin{array}{l}
J'\textrm{ is a sequence obtained from $J$ by deleting at least }\\
\textrm{one index and permuting the resulting sequence cyclicly}
\end{array}
\right\}.\]

Then we have the following theorem. 

\medskip
\begin{thm}\label{inv}
For any sequences $Ii$ and $J$ of length at most $q$, both $\Delta_{(C,\mathbf{p})}(J)$ and 
$\gcd\{m_{(C,\mathbf{p})}(Ii),\Delta_{(C,\mathbf{p})}(Ii)\}$ are invariants of 
$C$. In particular they are independent of the choice of 
$\mathbf{p}$. 
\end{thm}

\medskip
We define 
\[\nu_C(Ii):=\gcd\{m_{(C,\mathbf{p})}(Ii),\Delta_{(C,\mathbf{p})}(Ii)\}.\]
As for link diagrams, since the condition $s<q$ is not essential, 
by Theorem~\ref{inv}, we have the invariants $\nu_C(I)$ for all sequences $I$, and 
call them {\em Milnor invariants of the cut-diagram} 
$C$.

For a surface-link $L$, let $C$ be a cut-diagram obtained from a diagram of $L$.
For a sequence $I$, we define 
\[\nu_L(I):=\nu_C(I). \]
Then we have the following theorem. 

\medskip
\begin{thm}
If two surface-link diagrams of $L$ and $L'$ are deformed into each other by Roseman moves, 
then $\nu_L(I)=\nu_{L'}(I)$ for any sequence $I$. 
\end{thm}

\medskip
This theorem implies that $\nu_L(I)$ of all $I$ are invariants of $L$. 
We call the invariants {\em Milnor invariants of $L$}.

\medskip
\subsection{Based cut-diagrams and Chen map}

In this section, we explain how to calculate Milnor invariants of 2-dimensional cut-diagrams. 

Let $(C,\mathbf{p})~(C=(\Sigma,D,f))$ be a based cut-diagram. 
For each region $a_{ij}$ on the $i$th component $\Sigma_i$ of  $\Sigma$, 
we choose an oriented curve $\gamma_{ij}$ from $p_i$ to a point $p_{ij}$ in $a_{ij}$. 

While running the curve $\gamma_{ij}$ from $p_i$ to $p_{ij}$,
we obtain a word $v_{ij}$ by arranging the labels of arcs, powered by its sign, that appear in order, 
where the signs are defined in  Section~4.5. 

Let $F$ be the free group generated by $a_{i0}=\alpha_i~(i=1,...,n)$. 
For a positive integer $q$, by using $\{v_{ij}\}_{ij}$ above, we inductively  define a homomorphism 
\[\eta_{q}=\eta_{q}(C,\mathbf{p}):\widetilde{F}\longrightarrow F\] 
as follows.\footnote{The idea of this map was inspired by Chen's paper \cite{Chen}.}
\[(\mathrm{i})~\eta_{1}(a_{ij})=\alpha_{i};~~~
(\mathrm{ii})~\eta_{q+1}(a_{i0})=\alpha_{i},~\eta_{q+1}(a_{ij})=\eta_{q}(v_{ij}^{-1})\alpha_{i}\eta_{q}(v_{ij}) \hspace{1em} (1\leq j\leq r_i). \]

We call the map $\eta_q$ a {\em Chen map of $(C,\mathbf{p})$}. 
The map $\eta_q$ is different from the case of link diagrams, and  
it depends on the choices of not only base points but also the curves 
$\{\gamma_{ij}\}_{ij}$.

Lemma~\ref{Chen}~(1) holds for $\eta_{q}(C,\mathbf{p})$ as well. 
Moreover, the following theorem holds.

\medskip
\begin{thm}\label{Cut-CM}~
\begin{enumerate}
\item The Chen map $\eta_q(C,\mathbf{p})$ induces the following isomorphism. 
\[N_q(C,\mathbf{p})\cong\langle\alpha_1,\ldots,\alpha_n~|~ 
\{[\alpha_i,\eta_q(w(l_{ij}))]~|~i=1,...,n,~j=1,...,2g_i\},
\Gamma_qF\rangle.\]
\item For any loop $l$ with base point $p_i$, 
$\xi_q(w(l))=\eta_q(w(l)){W\Gamma_qF}$, 
where $W$ is the normal closure of $F$ that contains $\{[\alpha_i,\eta_q(w(l_{ij}))]~|~i=1,...,n,~j=1,...,2g_i\}$. 
\end{enumerate}
\end{thm}

\medskip
By Theorem~\ref{Cut-CM} (2), it can be seen that Milnor invariants of $C$ 
are obtained by replacing $\xi_q(w(l))$ with $\eta_q(w(l))$ in the definition.
Hence, $\eta_q$ gives Milnor invariants of cut-diagrams.

\section{Further generalization}

For any positive integer $m$, we define $m$-dimensional cut-diagrams, and 
their equivalence relation, {\em cut-concordance}, which corresponds to a 
generalization of the usual link concordance for $m$-dimensional links.

\medskip
\subsection{$m$-dimensional cut-diagrams and cut-concordance}

For a {\em diagram}\footnote{See \cite{Rose2}, \cite{Rose3}, or \cite{PR}.} $Y$
of compact $(m-1)$-manifold in a compact $m$-manifold $X$, we consider a map 
\[f:\{\text{~regions of $Y$}\}\longrightarrow\{\text{~regions of $X-Y$}\}\]
that satisfies \lq  certain gluing conditions\rq.\footnote{For the detailed definition, see \cite{AMY}.}
Here \lq regions of $Y$\rq\ correspond to arcs in diagrams in the case of $2$-dimensional cut-diagrams. 
We call the triple $(X,Y,f)$ an {\em $m$-dimensional cut-diagram}. 

Let $X$ be a closed $m$-manifold.\footnote{It is not necessary that $X$ is closed. 
But for simplicity we assume that.}
Two $m$-dimensional cut-diagrams $C_0=(X,Y_0,f_0)$ and $C_1=(X,Y_1,f_1)$ are 
{\em cut-concordant} if there is an $(m+1)$-dimensional cut-diagram $C=(X\times[0,1], Y,f)$
that satisfies the following two conditions for each $\varepsilon(=0,1)$:
  \begin{enumerate}
  \item There is an orientation preserving diffeomorphism 
    $\psi_\e:X\longrightarrow X\times\{\e\}$ such that
    $\psi_\e(Y_\e)=Y\cap\big(X\times\{\e\}\big)$.
   \item The following diagram commutes: 
\end{enumerate}
\[
\xymatrix{
\{\text{~regions of $Y_\e$}\} \ar[r]^{\hspace{-6ex}f_{\varepsilon}} 
\ar[d]^{\psi_\e^{Y_\e}} & ~~\{\text{~regions of $X-Y_\e$}\}~~~~~~~~ 
\ar[d]^{\psi_\e^{X-Y_\e}} \\
\{\text{~regions of $Y$}\} \ar[r]^{\hspace{-6ex}f} & \{\text{~regions of $(X\times [0,1])-Y$}\}. &
}
\] 
Here $\psi_\e^{Y_\e}$ (resp. $\psi_\e^{X-Y_\e}$) 
is a map induced by $\psi_\e$.
The condition~(1) implies that $Y_0\subset X$ and $Y_1\subset X$ cobound 
$Y\subset X\times[0,1]$, 
and the condition~(2) implies that 
the gluing on $C_0$ and $C_1$ can be extended to the gluing on $C$.

The cut-concordance is an equivalence relation on cut-diagrams. 
This is a generalization of {\em link concordance} on {\em $m$-dimensional links}, that 
are embedded $m$-manifolds in the $(m+2)$-space. 
In fact, the following proposition holds. 

\medskip
\begin{prop}\label{prop:ConcordanceSameSame}
If two $m$-dimensional links $L$ and $L'$ are link concordant, then two 
cut-diagrams obtained from link diagrams of $L$ and $L'$ are 
cut-concordant. 
\end{prop}

\medskip
\subsection{Milnor invariants for higher dimensional links.}

For an $m$-dimensional cut-diagram $C=(X,Y,f)~(m\geq2)$, 
the {\em Milnor invariants} $\nu_C(I)$ can be defined in the same way as for $m=2$.
When $m=1$, as mentioned in Subsection 4.4.1, 
the virtual link $L_C$ obtained from $C$ is uniquely determined as a welded link. 
Threfore we define $\omu_C(I):=\omu_{L_C}(I)$. 
Then we have the following theorem.

\medskip
\begin{thm}\label{cutconco}
If two $m$-dimensional cut-diagrams $C$ and $C'$ are cut-concordant, then 
for any sequence $I$, the following (1) and (2) hold.\\
\centerline{(1) $\omu_C(I)=\omu_{C'}(I)$ if $m=1$, ~~~~~~
(2) $\nu_C(I)=\nu_{C'}(I)$ if $m\geq 2$.}
\end{thm}

\medskip
Let $C$ be a cut-diagram obtained from a diagram of an $m$-dimensional link $L~(m\geq2)$. 
For any sequence $I$, we define 
\[\nu_L(I):=\nu_C(I).\]
Then by Proposition~\ref{prop:ConcordanceSameSame} and Theorem~\ref{cutconco}
we have the following corollary. 

\medskip
\begin{cor}
If two $m$-dimensional links $L$ and $L'$ are link concordant, then 
for any sequence $I$,  the following (1) and (2) hold.\\
\centerline{(1) $\omu_L(I)=\omu_{L'}(I)$ if $m=1$, ~~~~~~
(2) $\nu_L(I)=\nu_{L'}(I)$ if $m\geq 2$.}
\end{cor}

\medskip
This implies that for an $m$-dimensional link $L~(m\geq 2$), 
$\nu_L(I)$ is a link concordance invariant of $L$. 
We call $\nu_L(I)$ a {\em Milnor invariant of the $m$-dimensional link} $L$.


\end{document}